\documentclass[11pt]{article}
\usepackage[utf8]{inputenc}
\usepackage{amsbsy,amsfonts,amsmath,amssymb,amsthm
,bbm,enumerate,euscript,graphicx,color}

\usepackage{hyperref} 
\hypersetup{
    colorlinks=true,
    linkcolor=blue,
    citecolor=blue,
    urlcolor = magenta,
}

\def\R{\mathbb{R}}
\def\N{\mathbb{N}}
\def\cN{\mathbf{N}}
\def\Z{\mathbb{Z}}
\def\E{\mathbb{E}}

\def\V{\textbf{Var}}
\def\P{\mathbb{P}}

\def\I{\mathbf{1}}
\def\D{\mathrm{d}\,}

\newtheorem{theo}{Theorem}
\newtheorem{coro}{Corollary}

\newtheorem{prop}{Proposition}

\theoremstyle{definition}
\newtheorem{defn}{Definition}
\newtheorem*{proo}{Proof}
\newtheorem{exam}{Example}

\numberwithin{equation}{section}
\newtheorem*{rema}{Remark}

\title{\bf {Non-hyperuniformity of Gibbs point processes with short range interaction}}
\author{David Dereudre\thanks{University of Lille, david.dereudre@univ-lille.fr} 
\and Daniela Flimmel\thanks{Charles University and University of Lille,
daniela.flimmel@karlin.mff.cuni.cz}}
\date{\today}

\begin{document}
\maketitle

\begin{abstract}
We investigate the hyperuniformity of marked Gibbs point processes with weak dependencies among distant points whilst the interactions of close points are kept arbitrary. Some variants of stability and range assumptions are posed on the Papangelou intensity in order to prove that the resulting point process is not hyperuniform. The scope of our results covers many frequently used models including Gibbs point processes with a superstable, lower-regular, integrable pair potential as well as the Widom--Rowlinson model with random radii or Gibbs point processes with interactions based on Voronoi tessellation and nearest neighbour graph.
\end{abstract}

\vspace{0.2cm}


\vspace{0.5cm}

\section{Introduction}
Point processes form the main building block of stochastic geometry serving to construct a broad spectrum of geometric models used for analysing spatial data, for instance in material science, particle physics, telecommunication or biology. In these days, the theory of point processes is available under a very abstract topology setting (\cite{DVJ08, SW08}), yet in the current paper, we restrict ourselves to point configurations being randomly sprinkled in $\R^d$. This restriction is natural and still allows to build many random geometric structures via using some connections among points (e.g. random graphs or tessellations), by adding some additional mark to each point (a number, set, etc.) or by the combination of the above stated methods.
 
 A stationary (translation invariant) point process $\Gamma$ in $\R^d$ is called \textit{hyperuniform} if it  exhibits small density fluctuations. Formally, that is, if the number of points of the process in a bounded domain  $N_\Lambda:= N_\Lambda(\Gamma)$ fluctuates in a lower order then the volume of the set, i.e. 
\begin{equation}\label{hyperuniformity}
    \lim_{\Lambda \nearrow \R^d} \frac{\V(N_\Lambda(\Gamma))}{|\Lambda|}=0,
\end{equation}
where by $\Lambda \nearrow \R^d$ we understand any sequence of increasing sets $\Lambda_1 \subseteq \Lambda_2 \subseteq \ldots$ tending to the whole space $\R^d$. Historically, hyperuniformity is related to studying the compressibility of matter.
In statistical mechanics, compressibility is a measure of the relative volume change of a fluid or solid as a response to a pressure, which can be understood as the relative size of fluctuations in particle density.
Nonetheless, knowing the order of the variance is of interdisciplinary interest. In material science, the concept of hyperuniformity enables the characterization of naturally organized structures such as crystals and quasicrystals (\cite{BH22, OSS17}). The quantitative characterization of fluctuations in the number of particles has a long history also in statistical physics. Under certain constraints on the two-point correlation function (in the physics literature referred to as sum rules), the authors in \cite{MY80} (with an extension in \cite{L83}) showed for infinite classical systems of particles with long range interactions that the variance of $N_\Lambda$ should increase as the surface area. A rigorous result for the one-dimensional Coulomb systems is shown in \cite{K74}. 
In the theory of random matrices, one-dimensional point patterns associated with the eigenvalues have been characterized by their density fluctuations (e.g. \cite{M91}). The measurement of galaxy density fluctuations is a standard approach to study the structure of the Universe (see \cite{P93}).
The concept of hyperuniformity and density fluctuations is identified across many other areas of fundamental science like computer science, number theory or biological sciences. Yet, the theoretical understanding of such systems is still limited. The first attempt to rigorously handle the concept of hyperuniformity in the physics of matter  was established in the seminal paper \cite{TS03} and the subsequent papers of the authors. The current state of the art is summarized in the survey \cite{T18}.

It seems that it is a very fundamental and practical
question of great interest to determine whether certain point process is hyperuniform or not. It is not surprising that it became a new and fashionable topic for many researchers in stochastic geometry and related fields. Especially, the study of Coulomb and Riesz gases in the context of the theory of point processes comes in the front of mathematical and physics scientific interest. Some answers were given for $d=2,3$ in \cite{CH19, LS18}. By estimating the structure factor, the authors in \cite{HGB23, KLH22} provide tests of hyperuniformity based on point samples. Not just the proof of hyperuniformity, but also a negative answer is of great value in stochastic geometry and spatial statistics. Nondegeneracy of the asymptotic variance \eqref{hyperuniformity} is one of the key assumptions for geometric central limit theorems (as in \cite{BYY19, PY01, SY13}). 

The question of hyperuniformity is imminent for some special point processes, especially if the exact distribution of $N_{\Lambda}$ is tractable. It is the case of a stationary Poisson point process with intensity $\lambda>0$, which obviously is not hyperuniform. Here, the number of points in any set is independent of the outside configurations. Generally, it is unclear what happens if interactions between the points inside the set and outside are introduced into the model. A standard approach to generate point patterns with interactions among points is to consider a Gibbs modification of some underlying measure, typically a Poisson point process through energy. The simplest model consists of taking into account only interactions between couples of points, known as a Gibbs point process with pairwise interaction. Also, interactions among $k$-tuples of points and other more complicated types of interactions are widely studied in stochastic geometry and spatial statistics. In this paper, we show that short range Gibbs point processes are not hyperuniform. By a short range process, we simply mean a point process such that interactions among points weaken with the distance. This phenomenon can be interpreted through the Papangelou intensity. A point process $\Gamma$ has a Papangelou intensity $\lambda^*$ if for 
 any non-negative function $f$, we may write
 $$\E\left[\sum_{x \in \Gamma}f(x,\Gamma\setminus \{x\})\right]=\int \E f(x,\Gamma) \lambda^*(x,\Gamma)\D x.$$
Intuitively, we interpret $\lambda^*(x,\gamma)\D x$ as the conditional probability of observing a point in the infinitesimally small neighbourhood of $x$ given that $\Gamma$ agrees with a configuration $\gamma$ outside this neighbourhood. For example, if $\varphi:\R^d\times \R^d\to \R \cup \{+\infty\}$ measures interaction between two points, then a Gibbs point process with pairwise interaction has a Papangelou intensity of the form
$$\lambda^*(x,\gamma)= z e^{-\beta \sum_{y \in \gamma}\varphi(x,y)},$$
where $z>0, \beta \geq 0$ are usually called activity and inverse temperature. Now, a point process is short range if
\begin{equation}\label{ratio}
    \frac{\lambda^*(x,\gamma\cup \{y\})}{\lambda^*(x,\gamma)} \to 1
\end{equation}
fast enough as $\|x-y\| \to \infty$ for almost all configurations $\gamma$. The interpretation is straightforward. Any point $y$ in the configuration $\gamma$ plays only a negligible role when introducing a new point $x$ far away from $y$. For a Gibbs point process with pair potential $\varphi$, \eqref{ratio} simply translates to $\varphi(x,y)\to 0$ fast enough. 

 In this paper, we prove a kind of compressibility of any scale in the bulk of the interacting particle systems. Our main theorem is formulated for infinite-volume Gibbs point processes with a general Papangelou intensity $\lambda^*$. We claim that if the Papangelou intensity satisfies some moment conditions, it is enough to verify that
\begin{equation}\label{ShortRange}
    \int\left(\E \left|1-\frac{\lambda^*(0,\Gamma \cup \{y\})}{\lambda^*(0,\Gamma)}\right|^{\alpha}\right)^{1/\alpha} \D {y} <\infty
\end{equation}
for some $\alpha\geq 1$ to ensure that $\Gamma$ is not hyperuniform. 
 Specially for pairwise interaction with potential $\Phi$, \eqref{ShortRange} translates into the integrability condition
 $$\int |1-e^{-\beta \Phi(x)}|\D x<\infty,$$
 which, combined with superstability and lower-regularity as understood in \cite{R70} generates a non-hyperuniform point process. In the context of pair potentials, the authors in \cite{G67, R70} also provide some variance estimates, but only for the finite volume Gibbs measures.
A more general question than the initial one is addressed in \cite{XY15}, where the authors study volume order fluctuations of score functions (including the number of points). The setting is for short range Gibbs processes that are dominated by a Poisson point process. It corresponds to purely repulsive interactions, which can be translated to the Papangelou intensity by assuming $\lambda^*\geq C>0$ everywhere. Such a setting is also discussed here in Section \ref{SS_3.2}. For our main theorem, we consider more general interactions that can be arbitrary in a short range and vanishing in a long range. It covers also processes that cannot be coupled with a Poisson point process.  The techniques in both papers are completely different. To show our results, we use an approach based on the GNZ formalism instead of the algorithmic construction of the Gibbs point processes used in \cite{XY15}

The organization of the paper goes as follows. In Section \ref{S2}, we aim to introduce all the necessary notation covering the theory of point processes with marks and present two versions of our main result. First, a non-hyperuniformity result for Gibbs point processes with pairwise interactions or interactions that become deterministic for distant points. Next, we state a theorem for Gibbs point processes such that the interactions are random everywhere, but also weak at long distances, which is suitable for a large variety of geometrically based interactions. The following Section \ref{S3} provides some further investigation and a couple of hints to verify Assumption \ref{ShortRange}. It is followed by Section \ref{S4}, where we provide a comprehensive application of our main results; first discussing pairwise interaction and then more general geometrical interaction. Finally, in Section \ref{S5}, we give rigorous proofs of our main statements. 

\tableofcontents

\section{Main results}\label{S2}

\subsection{Notation}

To describe marked Gibbs point processes in the space $\R^d$, we proceed as follows. Let $\mathcal{B}$ be the Borel $\sigma-$algebra on $\R^d$ and $\lambda$ the Lebesgue measure on $(\R^d, \mathcal{B})$. Furthermore, we introduce a complete, separable mark space $\mathbb{M}$, equipped with the associated Borel $\sigma$-algebra $\mathcal{B}_{\mathbb{M}}$ and a finite measure $\lambda_{\mathbb{M}}$. 

We consider point configurations on the space $E:= \R^d \times \mathbb{M}$ with $\sigma$-field $\mathcal{E}:= \mathcal{B}\times \mathcal{B}_{\mathbb{M}}$, where each point in $\R^d$ is associated with a mark belonging to $\mathbb{M}$. The reference measure on $E$ is the product measure $\nu:= \lambda\otimes \lambda_{\mathbb{M}}$. The marks are mutually independent random variables whose distribution $\mathbb{Q}_{\mathbb{M}}$ does not depend on the location.

Let $\cN$ denote the set of all locally finite configurations in $E$, i.e.
$$\cN:= \left\lbrace \gamma \subset E; |\gamma \cap (\Lambda\times\mathbb{M})|<\infty, \, \forall \Lambda\subset \R^d \text{ bounded}\right\rbrace$$
Moreover, we denote by $\cN_f$ the subset of $\cN$ consisting of finite configurations. We endow $\cN$ with $\mathcal{N}$, which is the smallest $\sigma$-field such that all projections $\gamma\mapsto \gamma\cap B$ are measurable for all $B \in \mathcal{E}$. For a point configuration $\gamma \in \cN$ and a fixed set $\Lambda \subset \R^d$, we denote by $\gamma_\Lambda$ the restriction of $\gamma$ to the set $\Lambda\times \mathbb{M}$, i.e. $\gamma_\Lambda := \gamma \cap (\Lambda \times \mathbb{M})$. 

By a point process, we understand a probability measure $\Gamma$ on $(\cN,\mathcal{N})$. For $u \in \R^d$, we let $\tau_u:(x,m)\mapsto (x+u,m)$ be the shift in the position coordinate. If $\gamma=\{\mathbf{x}_1,\mathbf{x}_2,\ldots\}\in \cN$, we write $\tau_u \gamma = \{\tau_u \mathbf{x}_1,\tau_u\mathbf{x}_2,\ldots\}$.
Moreover, we call a point process stationary, if its distribution is invariant w.r.t. $\tau_u$, i.e $\Gamma \stackrel{D}{=}\tau_u\Gamma, \forall u \in \R^d$.

For a bounded set $\Lambda\subset \R^d$ and $\gamma \in \cN$, we denote the number of points of $\gamma$ occurring in $\Lambda$ by $N_\Lambda:=N_\Lambda (\gamma)=\sum_{\mathbf{x}\in \gamma}\I_{\Lambda\times\mathbb{M}}(\mathbf{x})$. If $\Gamma$ is a marked point process, then $N_\Lambda(\Gamma)$ is a random variable with values in $\N\cup\{0\}$.

\subsection{Gibbs point process}
\begin{defn}\label{Def_Gibbs}
    Let $\lambda^*:E\times\cN \to \R$ be some measurable function. We call $\Gamma$ a \textit{Gibbs point process} associated with the \textit{Papangelou intensity} $\lambda^*$ if for all positive measurable $f:E\times\cN\to \R$ it solves the Georgii--Nguyen--Zessin equations
\begin{equation}\tag{GNZ}\label{GNZ}
    \int_{\cN} \sum_{\mathbf{x}\in \gamma} f(\mathbf{x},\gamma\setminus \{\mathbf{x}\}) \Gamma(\D \gamma)=\int_{\cN} \int_{E} f(\mathbf{x},\gamma) \lambda^*(\mathbf{x},\gamma)\nu(\D \mathbf{x}) \Gamma(\D \gamma).
\end{equation}
\end{defn}

\begin{rema}[The form of the Papangelou intensity]
Usually, the Papangelou intensity is given in a form      \begin{equation}\label{Gibbs}
\lambda^*(\mathbf{x},\gamma) = z \exp\{-\beta h(\mathbf{x}, \gamma)\},
\end{equation}
where $h:E\times\cN \to \R$ is called local energy, $z>0$ is the activity parameter and $\beta\geq 0$ is the inverse temperature.

We only consider here stationary Gibbs point processes. In this case, the related Papangelou intensity $\lambda^*$ is necessarily translation invariant simultaneously w.r.t. both coordinates, meaning that
$$\lambda^*(\mathbf{x}, \gamma)= \lambda^*(\tau_u\mathbf{x}, \tau_u\gamma), \quad \forall u \in \R^d.$$ 
 \end{rema}

\begin{rema}
The existence of a point process from Definition \ref{Def_Gibbs} is not discussed in the present paper, it is always implicitly assumed that we are given a well-defined Gibbs point process. However, note that \eqref{ratio} is related to the condition for existence as in \cite{DV20}. Some relevant existence results are mentioned throughout the text. In addition, the point process $\Gamma$ may not be uniquely determined by \eqref{Gibbs}. If we make a claim about a Gibbs process with Papangelou intensity $\lambda^*$, we mean that it holds for all processes corresponding to this Papangelou intensity.
\end{rema}

\subsection{Main results}
Our first result concerns point processes whose are short range and, moreover, the ratio \eqref{ratio} has an upper bound which is uniform in $\cN$. An important class of short range Gibbs point processes is given by those determined by an integrable pair potential (see Section \ref{Ex_PP} for details). For multibody potentials, we replace the integrability condition by a suitable deterministic bound, which is integrable in the sense of \cite{G67, R70} except we do not force the potentials to be integrable around the origin.
\begin{theo}\label{thm1}
    Let $\Gamma$ be a stationary Gibbs point process on $\R^d$ with Papangelou intensity $\lambda^*$ and let the following assumptions are satisfied
    \begin{enumerate}
    \item[(a1)]
      For all $m \in \mathbb{M}$, $\E\lambda^*((0,m),\Gamma)^2<\infty$. 
    \item[(a2)]  Assume there exists a function $\phi:\R^d\to \R$ and $\delta\geq 0$ such that $\int_{\R^d\setminus B(0,\delta)}\phi(x)\D x<\infty$ and for any $\|x-y\|>\delta$,
\begin{equation*}
    \left|1-\frac{\lambda^*(\mathbf{x},\gamma\cup \{\mathbf{y}\})}{\lambda^*(\mathbf{x},\gamma)}\right|\leq\phi(x-y).
    \end{equation*}
    \end{enumerate}
    Then there is a constant $\mathcal{C}_{nhyp}>0$ not depending on $\Lambda$ such that 
    \begin{equation}\label{VarBound1}
    \frac{\V(N_\Lambda(\Gamma))}{|\Lambda|}\geq \mathcal{C}_{nhyp}>0.
\end{equation}
    In particular, $\Gamma$ is not hyperuniform.
\end{theo}

Theorem \ref{thm1} applies to point processes where the interaction among distant points is deterministically bounded. This, however, is often too restrictive for processes with random range of interaction (e.g. processes with energies based on Voronoi tessellation). For such processes, we provide another result.

\begin{theo}\label{thm2}
Let $\Gamma$ be a stationary Gibbs point process on $\R^d$ with Papangelou intensity $\lambda^*$.
Assume that there are some $\alpha_1, \alpha_2 > 1, \frac{1}{\alpha_1}+\frac{1}{\alpha_2} = 1$  and $\delta\geq 0$ such that 

\begin{enumerate}
    \item[(A1)]$  \int_{\mathbb{M}} \E |\lambda^*((0,m),\Gamma)|^{2\alpha_1}\lambda_{\mathbb{M}}(\D m) < +\infty$
    \item[(A2)] $\int_{\R^d\setminus B(0,\delta)}\left(\int_{\mathbb{M}^2}\E \left|1-\frac{\lambda^*((0,m_1),\Gamma \cup \{(y,m_2)\})}{\lambda^*((0,m_1),\Gamma)}\right|^{\alpha_2}\D\lambda_{\mathbb{M}}(m_1,m_2)\right)^{1/\alpha_2} \D y <\infty.$
\end{enumerate}
Then there exists $\mathcal{C}_{nhyp}>0$ not depending on $\Lambda$ such that
\begin{equation}\label{VarBound2}
    \frac{\V(N_\Lambda(\Gamma))}{|\Lambda|}\geq \mathcal{C}_{nhyp}.
\end{equation}
In particular, $\Gamma$ is not hyperuniform.
\end{theo}
The proofs of Theorems \ref{thm1} and \ref{thm2} are postponed to Section \ref{S5}. The constant $\mathcal{C}_{nhyp}$ can be given in a closed form directly from the proof. We provide some lower bounds for those constants in special cases in Section \ref{S4}.

\begin{rema}[Assumptions (A1), (A2)]
 Note that $\E\lambda^*(\mathbf{0},\Gamma) := \lambda$ is the intensity of $\Gamma$. Loosely speaking, Assumption (A1) of Theorem \ref{thm2} prevents the point process $\Gamma$ from having too many points in a unit window. In other words, it forces some stability assumptions. (e.g. superstability for pair potentials). On the other hand, Assumption (A2) states that the interaction among points of $\Gamma$ become negligible with the distance. Again, we do not force the interactions among close points to be bounded. A similar interpretation applies for Assumptions (a1), (a2) from Theorem \ref{thm1}.
\end{rema}

\begin{rema}[Shape of $\Lambda$]
Note that the proof and the lower bounds \eqref{VarBound1} \eqref{VarBound2} do not depend on the shape of the window $\Lambda$. 
\end{rema}

\begin{rema}[Unmarked case]
In the unmarked case, one can chose the mark space $\mathbb{M}$ with just one atom $m$ and set $\lambda_{\mathbb{M}}(m)=1$. The Papangelou intensity obviously does not depend on the choice of $m$, so we may set $\Gamma'= \{x_i; \,(x_i,m_i)\in \Gamma\}$ to be the unmarked point process and write $\lambda^*((x,m),\Gamma)=\lambda^*(x,\Gamma')$ a.s. Then the assumptions of Theorem \ref{thm2} turn to
\begin{enumerate}
    \item[(A1)]$ \E |\lambda^* (0, \Gamma')|^{2\alpha_1}< +\infty$,
    \item[(A2)] $\int_{\R^d\setminus B(0,\delta)}\left(\E \left|1-\frac{\lambda^*(0,\Gamma' \cup \{y\})}{\lambda^*(0,\Gamma')}\right|^{\alpha_2}\right)^{1/\alpha_2} \D {y} <\infty.$
\end{enumerate}
\end{rema}

\section{Conditions for (a1), (a2), resp. (A1), (A2)}\label{S3}
In this section, we provide a short cookbook of conditions posed on the Papangelou intensity $\lambda^*$ that guarantee the validity of Assumptions (a1), (a2) of Theorem \ref{thm1}, resp. Assumptions (A1), (A2) of Theorem \ref{thm2}.

\subsection{Stability and range of interaction}
\begin{defn}[Local stability]\label{D_LocStable}
The Papangelou intensity $\lambda^*$ is called
\begin{itemize}
    \item \textit{locally stable from above} if there is a constant $C_1 <\infty$ such that $\lambda^*(\mathbf{x},\gamma) \leq C_1$ for all $\mathbf{x}\in E$ and $\gamma \in \cN$,
     \item \textit{locally stable from below} if there is a constant $C_2> 0$ such that $\lambda^*(\mathbf{x},\gamma) \geq C_2$ for all $\mathbf{x}\in E$ and $\gamma \in \cN$,
     \item \textit{double locally stable} if it is simultaneously locally stable from above and from below.
\end{itemize}

\end{defn}
Note that a point process with locally stable from above Papangelou intensity satisfies trivially Assumption (a1) of Theorem \ref{thm1}, resp. (A1) of Theorem \ref{thm2}. Yet, it is one of the very frequent assumptions in the literature.

\begin{defn}[Range of interaction]\label{Def_range}
Let $\Gamma$ be a Gibbs point process with Papangelou intensity $\lambda^*$. Then $\Gamma$ has
\begin{itemize}
    \item \textit{finite range of interaction} if there exists $R>0$ such that for all $\gamma \in \mathcal{N}$
$$\lambda^*(\mathbf{x},\gamma)= \lambda^*(\mathbf{x}, \gamma \cap B(x,R)\times \mathbb{M}), \quad \forall \mathbf{x}:=(x,m)\in E.$$

\item  \textit{random finite range of interaction} if for all $\mathbf{x}:=(x,m) \in E$ there is an almost surely finite random variable $R_{\mathbf{x}}:= R_{\mathbf{x}}(\gamma)$ such that
$$\lambda^*(\mathbf{x},\gamma) = \lambda^*(\mathbf{x}, \gamma \cap B(x,R_{\mathbf{x}})\times\mathbb{M}), \quad \text{ for $\Gamma$-a.a. }\gamma.$$

\item \textit{decreasing range of interaction} if 
\begin{equation}\label{decreasingRange}
R_\mathbf{x}(\gamma) \geq R_\mathbf{x}(\gamma \cup \{\mathbf{y}\})\quad \forall \mathbf{x}, \mathbf{y} \in E, \gamma \in \cN.
\end{equation}  
\end{itemize}
\end{defn}

\begin{rema}
    Note that \eqref{decreasingRange} is in fact a natural condition arising from many geometric models including Voronoi tessellation, Delaunay triangulation or $k$-nearest neighbour graphs.
\end{rema}

\begin{prop}\label{L1}
Let $\Gamma$ be a Gibbs point process associated with Papangelou intensity $\lambda^*$.
\begin{enumerate}
    \item[(i)]  A Gibbs point process with locally stable from above Papangelou intensity satisfies Assumption $(a1)$ of Theorem \ref{thm1}, resp. $(A1)$ of Theorem \ref{thm2} with any $\alpha_1>0$.
    \item[(ii)] if $\Gamma$ has a random finite range of interaction and there exists a function $f:\R^+\to \R$ such that $\lambda^*(\mathbf{x},\gamma)\leq f(R_{\mathbf{x}}(\gamma))$ for all $\mathbf{x} \in E$ and $\Gamma$-a.a. $\gamma$. If $\E  f(R_\mathbf{0})^{2\alpha_1}<\infty$ then $\Gamma$ satisfies Assumption $(A1)$ of Theorem \ref{thm2}.
    \item[(iii)] 
    If $\Gamma$ has a finite range of interaction, then it satisfies Assumption (a2) of Theorem \ref{thm1} and Assumption (A2) of Theorem \ref{thm2} for any $\alpha_2>1$.
     \item[(iv)] If $\lambda^*$ is double locally stable and $\Gamma$ has a decreasing random finite range of interaction such that $\E R^\alpha_{\mathbf{0}} <\infty$ for some $\alpha> d$, then Assumption $(A2)$ of Theorem \ref{thm2} is satisfied for any $\alpha_2 \leq \alpha/d$. 
\end{enumerate}
    
\end{prop}
\begin{proo}
\begin{enumerate}
    \item[(i)] Since $\lambda^*$ is uniformly bounded, any moment is finite.
    \item[(ii)] For any $m \in \mathbb{M}$, we have
$$\E\lambda^*((0,m),\Gamma)^{2\alpha_1}\leq\E  f(R_{(0,m)})^{2\alpha_1}<\infty. $$
Since $\lambda_{\mathbb{M}}$ was assumed to be finite, (A1) is satisfied for any mark distribution.
    \item[(iii)] Consider $\delta = R$ and $\phi \equiv 0$ in Theorem \ref{thm1}, resp. $\delta = R$ in Theorem \ref{thm2}.
    \item[(iv)] 
    Take $\alpha_2 \in (1,\alpha/d)$. By using the Markov inequality, we arrive at
\begin{align*}
    \int_{\R^d}&\left(\int_{\mathbb{M}^2}\E \left|1-\frac{\lambda^*((0,m_1),\Gamma \cup \{(y,m_2)\})}{\lambda^*((0,m_1),\Gamma)}\right|^{\alpha_2}\D \lambda^2_{\mathbb{M}}(m_1,m_2) \right)^{1/\alpha_2} \D y\\ 
    &\leq  \left(1+\frac{C_2}{C_1}\right)\int_{\R^d} \left[ \int_{\mathbb{M}^2}\P \left(|y|\leq R_{(0,m_1)}(\Gamma\cup\{(y,m_2)\})\right)\D \lambda^2_{\mathbb{M}}(m_1,m_2)\right]^{1/\alpha_2}\D y\\
    &\leq  \left(1+\frac{C_2}{C_1}\right)\lambda_{\mathbb{M}}(\mathbb{M})^{1/\alpha_2}\int_{\R^d} \left[ \int_{\mathbb{M}}\P \left(|y|\leq R_{(0,m_1)}(\Gamma)\right)\D \lambda_{\mathbb{M}}(m_1)\right]^{1/\alpha_2}\D y\\
    &\leq  \left(1+\frac{C_2}{C_1}\right)\lambda_{\mathbb{M}}(\mathbb{M})^{1/\alpha_2}\int_{\R^d} \frac{1}{|y|^{\alpha/\alpha_2}} \D y \left( \int_{\mathbb{M}}\E R_{(0,m_1)}^{\alpha}(\Gamma)\D \lambda_{\mathbb{M}}(m_1)\right)^{1/\alpha_2}<\infty.
\end{align*}
\end{enumerate}

\qed
\end{proo}

In conclusion, the task of verifying assumptions of Theorems \ref{thm1} and \ref{thm2} often translates to estimating the moments, resp. the tail probabilities of the radius of interactions. In more than a few situations, one could benefit from the literature on stabilization and stochastic relation to the Poisson point process or other random structures. The latter is described in the following section.

\subsection{Stochastic comparison tools}\label{SS_3.2}
In this section, we explore tools to estimate the moments $\E R_{\mathbf{0}}(\Gamma)^\alpha, \alpha>0$ of the range of interaction needed in the previous section. Usually, for an infinite volume Gibbs point process, it is not a straightforward task, since we do not possess the local distribution of the number of points. We use stochastic comparison with some random object which is easier to handle.

For this reason, we consider the usual order on $\cN$. For $\gamma_1, \gamma_2 \in \cN$, we write $\gamma_1\leq \gamma_2$ if $\gamma_1(B)\leq \gamma_2(B)$ for all $B \in \mathcal{E}$. In the language of point sets, this means that $\gamma_1$ has fewer points than $\gamma_2$. A function $f:\cN\to \R$ is called \textit{increasing} if $f(\gamma_1)\leq f(\gamma_2)$ whenever $\gamma_1\leq \gamma_2$. Inversely, it is \textit{decreasing}, if $-f$ is increasing.

\begin{defn}
    We say that a point process $\Gamma_1$ is \textit{stochastically dominated} by a point process $\Gamma_2$ (denoted by $\Gamma_1\ll\Gamma_2)$ if $\int h \D \Gamma_1 \leq \int h \D \Gamma_2$ for all increasing functions $h$. Vise versa, we say that $\Gamma_1$ \textit{stochastically minorates} $\Gamma_2$.
\end{defn}
By the famous Strassen theorem, $\Gamma_1\ll\Gamma_2$ if and only if there is a coupling of $\Gamma_1$ and $\Gamma_2$ that is supported in the set $\{(\gamma_1,\gamma_2), \gamma_1\leq \gamma_2\}$. 

First, we recall the “Poisson sandwich inequality” from \cite{GK97} which stochastically connects a Gibbs point process with a stationary Poisson point process.

\begin{prop}\label{L_PoissonSandwich}
Let $\Gamma$ be a Gibbs point process associated with Papangelou intensity $\lambda^*$ and let $\Pi_\rho$ denote a stationary Poisson point process with intensity $\rho$.
\begin{itemize}
    \item[(i)] If $\lambda^*$ is locally stable from below, then $\Pi_{{C_1}}\ll\Gamma$,
    \item[(ii)] if $\lambda^*$ is locally stable from above, then $\Gamma\ll\Pi_{{C_2}}$,
\end{itemize}
where $C_1, C_2$ are the constants from Definition \ref{D_LocStable}.
\end{prop}

\begin{coro}\label{C01}
Let $\Gamma$ be associated with the Papangelou intensity $\lambda^*$ and assume that the corresponding range of interaction $R_{\mathbf{x}}$ is decreasing. If $\lambda^*$ is locally stable from below, then 
$$\P(R_\mathbf{x}(\Gamma)>r) \leq \P(R_\mathbf{x}(\Pi_{C_1})>r), \quad \forall r \geq 0,$$
where $\Pi_{C_1}$ is as in Proposition \ref{L_PoissonSandwich}.
\end{coro}

\begin{rema}[Increasing range of interaction]
 If the range of interaction is increasing with respect to the order on $\cN$, then Corollary \ref{C01} provides the inequality in the opposite direction. In this situation, the local stability of $\lambda^*$ from above produces similar upper bounds for the tail probabilities of $R_\mathbf{x}(\Gamma)$.
\end{rema}

 Consequently, for $\Gamma$ as in Corollary \ref{C01} and $f:\R \to \R$ an increasing function, we may verify the assumptions of Proposition \ref{L1}, (ii), and Proposition \ref{L1}, (iv), by computing
\begin{equation}\label{EstMoments1}
    \E f(R_{\mathbf{0}}(\Gamma))=\int_{0}^{\infty}\P(f(R_{\mathbf{0}}(\Gamma))>r)\D r \leq \int_{0}^{\infty}\P(f(R_{\mathbf{0}}(\Pi_{}))>r)\D r=\E f(R_{\mathbf{0}}(\Pi_{})).
\end{equation}
if $f:\R \to \R$ is some increasing function

Theorems \ref{thm1} and \ref{thm2} provide non-hyperuniformity result for Gibbs processes with interactions that are weak at long distances. As mentioned before, we assume nothing about the interactions among close points. Such points can generate an unboundedly large amount of energy (imagine Coulomb interaction) and, therefore, the value $\lambda^*(\mathbf{x},\gamma)$ can be arbitrarily close to zero, when $\mathbf{x}$ is at a small distance of $\gamma$, i.e. $d(\mathbf{x},\gamma):= \inf\{d(\mathbf{x},\mathbf{y}); \mathbf{y} \in \gamma\}$ is small, where $d(\mathbf{x},\mathbf{y})$ is some distance of two marked points. This situation does not allow one to use minoration by a Poisson point process. However, in this situation, we are able to construct a coupling with a Bernoulli field.

For this purpose, we introduce the mapping $I^{s}: \cN \to \{0,1\}^{\Z^d}$ as follows. We split the space $\R^d$ into a collection of disjoint cubes $\mathcal{D}:=\{\mathcal{D}_i; i\in \Z^d\}$ of a common side length $s>0$ such that there exists $k \in \Z^d$ with the origin being one of the vertices of $\mathcal{D}_k$. Then, we define 
$$I^s:\gamma \mapsto (\I\{N_{\mathcal{D}_k}(\gamma)\geq 1\})_{k \in \Z^d}, \quad \gamma \in \cN.$$
In fact, in many situations, we do not need to know the exact positions of the points of the process to estimate the range of interaction. The only necessary information is often that there is at least one point in a given region. That is exactly the meaning of $I^{\delta}$.

In order to formulate our next result, we need some partial order on $\{0,1\}^{Z^d}$. We write that $l_1 \leq l_2$ if $l_1,l_2 \in \{0,1\}^{\Z^d}$ are such that $l_{1,i}=1$ implies $l_{2,i}=1$ for all $i \in \Z^d$. As usual, we say that $R:\{0,1\}^{Z^d}\to \R$ is increasing if $R(l_1)\leq R(l_2)$ whenever $l_1\leq l_2$ and decreasing if $-R$ is increasing. We write $X\ll Y$ for two random variables $X,Y$ with values in $\{0,1\}^{\otimes \Z^d}$ if $\P(X_{i}=1; i \in J)\leq \P(Y_{i}=1, i \in J)$ for any $J \subset \Z^d$.

\begin{prop}\label{L4}
Assume there are constants $C_1, C_2, \delta>0$ such that $\lambda^*(\mathbf{x},\gamma)\leq C_1$ everywhere and $\lambda^*(\mathbf{x},\gamma)\geq C_2$ whenever $d(\mathbf{x},\gamma)\geq \delta$. If $\Gamma$ is the corresponding Gibbs point process, then for any $\varepsilon>0$ there exists $p>0$ such that
\begin{equation}\label{dom_Bernoulli}
  B < < I^{2\delta+\varepsilon}(\Gamma),  
\end{equation}
where $B$ is a random variable with values in $\{0,1\}^{\otimes \Z^d}$ such that $B_i, i \in \Z^d$ are independent with Bernoulli distribution $B(p)$.
\end{prop}
See Section \ref{Sec_P5} for the proof. 
\begin{rema}[On the constant $\varepsilon$]
    The constant $\varepsilon$ in Proposition \ref{L4} is rather just auxiliary for the proof. We need the side length of $\mathcal{D}_k$ to be slightly bigger then $2\delta$ so that we can fit another box $\mathcal{C}_k$ of side length $\varepsilon$ inside $\mathcal{D}_k$ such that it has distance exactly $\delta$ from the complement of $\mathcal{D}_k$. Depending on $\varepsilon$, though, one can optimize the value of $p$ for further application.
\end{rema}

\begin{coro}\label{C02}
    Let the range of interaction $R_{\mathbf{x}}$ be decreasing with respect to the order in $\cN$ and, moreover, there exists a decreasing function $R':\{0,1\}^{Z^d}\to \R$ such that
    $R_{\mathbf{x}}(\gamma)\leq R'(I^{\delta}(\gamma))$ for all $\gamma \in \cN$ and some $\delta>0$. Then, under the assumption and notation of Proposition \ref{L4}, we have
 $$\P(R_{\mathbf{x}}(\Gamma) > r)\leq \P(R'(B)> r) \quad \text{for all }r>0.$$
\end{coro}
Note that $R_{\mathbf{x}}$ and $R'$ are defined on different spaces, and so the assumption of being decreasing has a slightly different meaning.

We know exactly the distribution of $B$. Therefore, it is usually a simple task to compute $\P(R'(B)> r)$. Then, using Corollary \ref{C02}, we might estimate the moments of $R_{\mathbf{0}}(\Gamma)$ similarly as in \eqref{EstMoments1} for any increasing function $f:\R\to \R$ by
$$\E f(R_{\mathbf{0}}(\Gamma))\leq \E f(R'(B)).$$

Some application of both Corollary \ref{C01} and Corollary \ref{C02} shall be demonstrated in Section \ref{SS_Voronoi} and \ref{SS_NNG}.

\section{Examples}\label{S4}

\subsection{Pair potentials}\label{Ex_PP}
A classical example of a Gibbs point processes in $\R^d$ is the model with pairwise interactions. In this section, we omit the marks.
\begin{defn}[Pair potential]
A Gibbs point process has a \textit{pair potential}, if 
there is a measurable, symmetric function $\Phi: \R^d \to \R\cup\{+ \infty\}$ such that the Papangelou intensity has a form
\begin{equation*}\label{PP}
   \lambda^*(x,\gamma)= ze^{-\beta \sum_{y \in \gamma} \Phi(x-y)}, \quad  x \in \R^d, \gamma \in \cN. 
\end{equation*}
\end{defn}

In the following definition, we recall a classical stability assumption in the pairwise model originating in \cite{R70}. Note that these conditions altogether guarantee existence of a Gibbs point process with pair potential $\Phi$.

 \begin{defn}[Superstable pairwise interaction]
    Let $\Gamma$ be a Gibbs point process with pair potential $\Phi$. For any finite configuration $\gamma \in \cN_f$, we define the energy of this configuration by
    $$H(\gamma):= \sum_{\substack{\{x,y\} \subset \gamma\\ x \neq y}}\Phi(x-y).$$
    We say that $\Phi$ is
    \begin{itemize}
        \item \textit{superstable}, if for any $\Lambda \subset \R^d$ bounded there exist constants $A>0, B\geq 0$ such that 
        $$H(\gamma_\Lambda)\geq A N_{\Lambda}(\gamma)^2 - B N_\Lambda(\gamma),\quad \forall \gamma \in \cN,$$
        \item \textit{lower regular}, if there is a positive, decreasing function $\varphi: [0,\infty)\to \R$ such that
        $$\int_0^{\infty} x^{d-1}\varphi(x) \D x <\infty$$
        and $\Phi(x)\geq - \varphi(|x|)$ for all $x \in \R^d$,
        \item \textit{integrable}, if
\begin{equation}\label{integrability}
    \int_{\R^d}\left|1-e^{-\beta\Phi(y)}\right|\D y <\infty.
\end{equation}
    \end{itemize}
 \end{defn}

\begin{rema}[Assumptions (a1), (a2) for superstable interactions]
       A Gibbs point process with an integrable pair potential $\Phi$ satisfies automatically Assumption (a2) of Theorem \ref{thm1} with $\phi=\Phi$ since
    $$\frac{\lambda^*(0,\gamma\cup \{y\})}{\lambda^*(0,\gamma)}=1-e^{-\Phi(y)}, \quad \forall \gamma \in \cN.$$
    Therefore, Assumption (a2) of Theorem \ref{thm1} reduces to the standard integrability condition of the pair potential $\Phi$. On the other hand, the existence of the second moment stems from the assumption of superstability and lower regularity (see Corollary 5.3 in \cite{R70}).
    \end{rema}

\begin{coro}\label{C2}
A Gibbs point process with superstable, lower-regular and integrable pair potential is not hyperuniform.
\end{coro}

\begin{rema}[Bound on the asymptotic variance]
    Following the proof of Theorem \ref{thm1}, one can derive the lower bound for the asymptotic variance of Theorem \ref{thm1} for the Gibbs point process with integrable pair potential $\Phi$:
    \begin{equation}\label{PP_bound}
        \mathcal{C}_{nhyp}\geq \frac{(\E\lambda^*(0,\Gamma))^2}{\E\lambda^*(0,\Gamma)+
       \E\lambda^*(0,\Gamma)^2\int_{\R^d}\left|1-e^{-\beta\Phi(x)}\right|\D x }.
    \end{equation} 
\end{rema}

\begin{exam}
The class of Gibbs point processes with superstable, integrable pairwise interactions is large and covers a lot of standard examples, for instance
\begin{enumerate}
    \item  \textit{Strauss process}, i.e. a Gibbs point process $\Gamma$ with pair potential 
$$
\Phi(x-y)=
\begin{cases}
1, \quad \text{ if } |x-y|\leq R,\\
0, \quad \text{ otherwise}
\end{cases}
$$
for some $R\in [0,\infty)$. Let $\lambda:=\E\lambda^*(0,\Gamma)$ be the intensity of the process, then directly from \eqref{PP_bound}, we get the lower bound 
$$\mathcal{C}_{nhyp}\geq \frac{\lambda^2}{z+z^2|B(0,R)|(1-e^{-\beta})}.$$
  
  \item \textit{Riesz gases} with $s>d$, i.e processes with a pair potential of the form $\Phi(x)=\|x\|^{-s}, x \in \R^d, s \in \R$. The case $s\leq d$ (also Coulomb gas $s=d-2$) determines a~non-integrable pair potential.  

\item \textit{Lennard--Jones} pair potential given by 
$$\Phi(x)= A\|x\|^{-\alpha_1}- B\|x\|^{-\alpha_2}$$ 
for some $A, B >0$ and $\alpha_1> \alpha_2> d$.

\end{enumerate}
\end{exam}

\subsection{Widom--Rowlinson models}
\label{ex_WR}
First, we start with a simple model where the size of the balls is the same deterministic constant for all the points. Thus, for $R>0$, we define
$$L_R(\gamma):= \bigcup_{x \in \gamma} B(x,R), \quad \gamma \in \cN.$$
The Widom--Rowlinson point process is a Gibbs point process $\Gamma$ with an energy function
$$H(\gamma)= |L_R(\gamma)|,\quad \gamma \in \cN_f.$$
Then, the Papangelou intensity has a form
$$\lambda^*(x,\gamma)= z\exp\{-\beta(|L_R(\gamma\cup \{x\})|-|L_R(\gamma)|)\},\quad x \in \R^d, \gamma \in \cN$$
for some $z>0, \beta\geq 0$.
Clearly, $z e^{-\beta |B(0,R)|}\leq \lambda^*(x,\gamma)\leq z$, hence $\lambda^*$ possess all moments finite. Moreover, $\Gamma$ has a finite range of interaction that produces non-hyperuniformity by Proposition \ref{L1}, (i) and (iii). Alternatively, we can directly use Theorem \ref{thm1} with $\delta=0$ and $\phi(x)=I_{B(0,2R)}(x) e^{\beta |B(0,R)|}$. Similarly, a Gibbs point process with quermass interaction among balls with a deterministic size is non-hyperuniform, i.e. a Gibbs point process with energy
$$H(\gamma) = |L_R(\gamma)| + \text{Per}(L_R(\gamma))+ \chi(L_R(\gamma)),\quad \gamma \in \mathbf{N}_f,$$
where $\text{Per}$ is the perimeter and $\chi$ the Euler--Poincaré characteristic (the number of connected components minus the number of holes).

 Moreover, by using the above mentioned estimates for the Papangelou intensity and by following the steps in the proof of Theorem \ref{thm1}, we provide lower bound for the asymptotic variance
$$\mathcal{C}_{nhyp}\geq \frac{e^{-\beta|B(0,R)|}}{1+ze^{\beta|B(0,R)|}|B(0,2R)|}.$$

Now assume each point $x\in \Gamma$ is equipped with a non-negative random variable $R_x$ independently and with the same distribution $\mathbb{Q}$. Now, we are in the setting of marked Gibbs point process with $\mathbb{M}= \R_+$ and $\lambda_{\mathbb{M}}=\mathbb{Q}$. As in the previous example, we define 
$$L(\gamma):= \bigcup_{(x,R_x) \in \gamma} B(x,R_x), \quad \gamma \in \cN$$
and 
\begin{equation}\label{eq_WRR}
\lambda^*(\mathbf{x},\gamma)= z \exp\{-\beta(|L(\gamma\cup \{\mathbf{x}\})|-|L(\gamma)|)\}.
\end{equation}
Note that a Gibbs point process associated with Papangelou intensity \eqref{eq_WRR} is well defined. The existence of the process with interactions involving also other Minkowski functionals is proved in \cite{D09}.

\begin{coro}\label{C3}
Assume that $\E_{\mathbb{Q}} e^{\alpha \beta |B(0,1)| R^d}<\infty$ for some $\alpha>1$, then a Gibbs point process defined by the Papangelou intensity \eqref{eq_WRR} is non-hyperuniform.
\end{coro}
\begin{proo}
As in Example \ref{ex_WR}, $\lambda^*\leq z$. Hence, Assumption $(A1)$ of Theorem \ref{thm2} is trivially satisfied by Proposition \ref{L1}, (ii). It remains to show that Assumption (A2) is true for some $\alpha_2>1$. Furthermore, for any $\mathbf{x}:=(x,R_x), \mathbf{y}:=(y,R_y) \in E:=\R^d\times\R_+$, we have
$$\frac{\lambda^*(\mathbf{x},\gamma\cup\{\mathbf{y}\})}{\lambda^*(\mathbf{x},\gamma)}\leq e^{\beta|B(x,R_x)\cap B(y,R_y)|},$$
which does not depend on the configuration $\gamma$. Therefore, for any $\alpha_2>1$,
\begin{equation}\label{eq4.4}
\E \left|1-\frac{\lambda^*(\mathbf{x},\Gamma \cup \{\mathbf{y}\})}{\lambda^*(\mathbf{x},\Gamma)}\right|^{\alpha_2}\leq\left|1-e^{\beta|B(x,R_x)\cap B(y,R_y)|}\right|^{\alpha_2}=\left(e^{\beta|B(x,R_x)\cap B(y,R_y)|}-1\right)^{\alpha_2}.
\end{equation}
 Now, take $\zeta_1, \zeta_2>1$ such that $\zeta_1<\alpha$ and $1/\zeta_1+1/\zeta_2=1$. Further, set $\alpha_2 = \alpha/\zeta_1$ and take $\zeta>\alpha_2 \zeta_2 d$. We use successively \eqref{eq4.4}, H\"older inequality with $\zeta_1,\zeta_2$ and Markov inequality with $\zeta$ to get
\begin{align*}
     \int_{\R^d}&\left(\int_{\R_+\times\R_+} \E \left|1-\frac{\lambda^*((0,R_0), \Gamma\cup \{(y,R_y)\})}{\lambda^*((0,R_0),\Gamma)}\right|^{\alpha_2}\mathbb{Q}(\D R_0)\mathbb{Q}(\D R_y)\right)^{1/\alpha_2} \D y\\
     &\leq \int_{\R^d}\left(\int_{\R_+\times\R_+}  \left(e^{\beta|B(0,R_0)\cap B(y,R_y)|}-1\right)^{\alpha_2}\I\{ R_0+R_y\geq |y|\}\mathbb{Q}(\D R_0)\mathbb{Q}(\D R_y)\right)^{1/\alpha_2} \D y\\
     & \leq \int_{\R^d}\left(\E e^{\beta\alpha_2\zeta_1|B(0,R_0)|}\right)^{1/\alpha_2\zeta_1}\left(\P( R_0+R_y\geq |y|)\right)^{1/\alpha_2\zeta_2} \D y\\
     & =\left(\E e^{\beta\alpha|B(0,1)|R^d_0}\right)^{1/\alpha}\int_{\R^d}\left(\P( R_0+R_y\geq |y|)\right)^{1/\alpha_2\zeta_2} \D y\\
     &\leq \left(\E e^{\beta\alpha|B(0,1)|R^d_0}\right)^{1/\alpha}\int_{\R^d}\frac{\left(\E(R_0+R_y)^\zeta\right)^{1/\alpha_2\zeta_2}}{|y|^{\zeta/\alpha_2 \zeta_3}} \D y <\infty.
\end{align*}

\qed
\end{proo}

\begin{rema}[Gibbs particle process]
Generally, we may replace balls with random radii by random compact sets to define the Gibbs particle process. Assuming that diameters of these sets have suitable exponential moments, the resulting process is again non-hyperuniform. Consequently, a crucial assumption of positive asymptotic variance in \cite{BHLV20} is satisfied in order to study the limit behaviour of $U$-statistics of Gibbs particle processes.
\end{rema}

\subsection{Voronoi interaction}\label{SS_Voronoi}
In the previous example, we had direct control over the sizes of the discs through marks and hence over the range of interaction. Here, we may encounter cells that are very large, interfere with distant neighbourhood, and its distribution is generally unknown. To define a Gibbs point process with interactions among Voronoi cells, we denote for nonempty $\gamma\in \cN$ and $x \in \gamma$ the cell around $x$ by
$$C(x, \gamma):= \{z \in \R^d: \|z-x\|\leq\|z-y\| \text{ for all } y \in \gamma\setminus\{x\}\}.$$
The set $C(x, \gamma)$\index{$C(x, \gamma)$} represents those points in $\R^d$ such that $x$ is their nearest point among points in $\gamma$.
If $C(x,\gamma \cup\{x,y\})
\cap C(y,\gamma\cup\{x,y\})\neq \emptyset$, we write $x\stackrel{\gamma}{\sim} y$. At last, if $\gamma\in \cN$ and $x\notin \gamma$, we write $C(x,\gamma)$ for $C(x,\gamma\cup \{x\}).$ The set $C(x, \gamma)$ is in fact a closed convex set, since it can be written as an intersection of closed half-spaces. The set of all closed convex sets is denoted by $\mathcal{C}^d$.

We consider a function $\Phi:\mathcal{C}^d\to \R\cup \{\infty\}$ and write $\Phi(x,\gamma)$ for $\Phi(C(x,\gamma))$. For a finite configuration $\gamma \in \cN$, we consider energy
\begin{equation}\label{VoroEnergy}
    H(\gamma)= \sum_{x \in \gamma}\Phi(x,\gamma)\I\{|C(x,\gamma)|<\infty\}.
\end{equation}

The Gibbs point process associated to $H$ is a point process satisfying \eqref{GNZ} with the Papangelou intensity of the form
\begin{align}\label{VoroGPP}
\lambda^*(x,\gamma) :&= z \exp\left\lbrace-\beta \Phi(x,\gamma) -\beta \sum_{y \in \gamma}\left[\Phi(y,\gamma\cup \{x\})-\Phi(y,\gamma)\right]\right\rbrace.
\end{align}
In order to analyze Assumptions (A1) and (A2) of Theorem \ref{thm2}, we provide a list of assumptions posed on the function $\Phi$.

\begin{defn}
We say that the function $\Phi$ is
\begin{enumerate}
    \item  \textit{sub-additive} if $\Phi(C)\leq \Phi(C_1)+\Phi(C_2)$ whenever $C =C_1\cup C_2$, where $C,C_1,C_2 \in \mathcal{C}^d$,
    \item \textit{increasing} if $\Phi(C)\leq \Phi(C')$ whenever $C,C' \in \mathcal{C}^d$ and $C \subseteq C'$,
\item \textit{controlled by the volume} if there exists a constant $K>0$ such that $|\Phi(C)|\leq \min\{|C|,K\}$ for all $C \in \mathcal{C}^d$.
\end{enumerate}
\end{defn}
Note that a Gibbs point process associated with $\lambda^*$ in \eqref{VoroGPP} with $\Phi$ as above exists according to \cite{DGD12}. 
The conditions stated above on $\Phi$ allow us to use stochastic minoration by the Poisson point process in order to gain control over the size of the typical cell.

\begin{prop}\label{L7}
If $\Phi$ is sub-additive, increasing and controlled by the volume, then 
\begin{equation}\label{localEnergyBound2}
    z e^{-\beta K}\leq \lambda^*(x,\gamma)\leq z e^{\beta K} e^{\beta|C(x,\gamma)|}.
\end{equation}
\end{prop}

\begin{proo}
For any $y \in \gamma$, we have $C(y,\gamma \cup \{x\})\subseteq C(y,\gamma)$ while both sides are elements of $\mathcal{C}^d$ and hence $\Phi(y,\gamma \cup \{x\}) \leq \Phi(y,\gamma)$ since $\Phi$ is increasing. Consequently, 
$$ \sum_{y \in \gamma}\left[\Phi(y,\gamma\cup \{x\})-\Phi(y,\gamma)\right]\leq 0.$$
Moreover, $\Phi$ is uniformly bounded from above by $K$ yielding $-\beta \Phi(x,\gamma)\geq -\beta K$. Therefore, the lower bound is proved.

To show the upper bound, we define for $\gamma \in \cN$ and $y\stackrel{\gamma}{\sim} x$ the set $K_y:= C(y,\gamma)\setminus C(y,\gamma\cup\{x\})\in \mathcal{C}^d.$ 
Then
$$
    C(x,\gamma)= \bigcup_{y\stackrel{\gamma}{\sim} x} K_y,
$$  
   but also
\begin{equation}\label{VoroCellDivision}
\sum_{y\stackrel{\gamma}{\sim} x}|K_y|=|C(x,\gamma)|,\end{equation}  
 since $|K_y\cap K_{y'}|=0$ for any $y\neq y'$.

Note that it is enough to consider only the neighbouring points of $x$ in \eqref{VoroGPP}, that is, to sum over $\{y \in \gamma; x\stackrel{\gamma}{\sim} y\}$. Otherwise, adding (or removing) a point $x$ into the configuration $\gamma$ does not affect the shape of $C(y,\gamma)$, hence the value of $\Phi(y,\gamma)$. Using the sub-additivity of $\Phi$, control by the volume and \eqref{VoroCellDivision}, we arrive at
\begin{align*}
    \sum_{y \in \gamma}\left[\Phi(y,\gamma\cup \{x\})-\Phi(y,\gamma)\right]& = \sum_{y \in \gamma; x\stackrel{\gamma}{\sim} y}\left[\Phi(y,\gamma\cup \{x\})-\Phi(y,\gamma)\right]\\
    & = \sum_{y \in \gamma; x\stackrel{\gamma}{\sim} y}\left[\Phi(y,\gamma\cup\{x\})-\Phi(C(y,\gamma\cup\{x\})\cup K_y)\right]\\
    &\geq - \sum_{y \in \gamma; x\stackrel{\gamma}{\sim} y} \Phi(K_y)\\
    &\geq - \sum_{y \in \gamma; x\stackrel{\gamma}{\sim} y} |K_y| = - |C(x,\gamma)|.
\end{align*}
The latter combined with $-\Phi\leq  K$ produce the desired upper bound.

 \qed
\end{proo}

\begin{coro}\label{C4}
 Let $\Phi$ be some sub-additive, increasing and controlled by the volume function on $\mathcal{C}^d$. Then a Gibbs point process $\Gamma$ associated with the Papangelou intensity \eqref{VoroGPP} is non-hyperuniform for any $\beta < \beta_c$, where $\beta_c$ is the unique solution of the equation
 $$ \beta_c=z e^{-\beta_c K}C_d,$$
 where
 $$C_d:=\frac{1}{3}\frac{|B(0,1)|_{d-1}}{|B(0,1)|_d}\left(\sin^{d-1}(\pi/12)\cos(\pi/12)\frac{1}{d}+\int_{0}^{\pi/12}\sin^d(\theta)\D \theta\right).$$
\end{coro}
The proof is postponed to Section \ref{S5}.

\begin{rema}[Values of $C_d$ and $\beta_c$]
    The constant $C_d$ can be easily evaluated. It can be seen that it decreases with dimension $d$, for instance $C_1 = \frac{1}{6}$, $C_2= \frac{1}{36}$, $C_3\approx 0.006$, etc. As a consequence, the exact value of $\beta_c$ can be also given depending on $z$ and $K$. If $d=2$ and $z=K=1$, then $\beta_c \approx 0.03$.
\end{rema}

 \begin{rema}[General tessellations]
Similarly, we can formulate Corollary \ref{C4} for interactions based on more general tessellations such as the Laguerre and Johnson--Mehl types. In this context, the shape of the cells is determined additionally by marks, which justifies considering the marks in Theorem \ref{thm1} and Theorem \ref{thm2}. 

 \end{rema}

\subsection{Interactions based on $k$-nearest neighbour graph}\label{SS_NNG}
    For $k\in \N$ and $\gamma \in \cN$, denote by $v^i(x,\gamma)$ the $i$-th nearest neighbour of $x$ in $\gamma$, $i=1,\ldots,k$ and let $V^k(x,\gamma):=\{v^i(x,\gamma);i=1,\ldots,\min\{k,N(\gamma)-1\}\}$ be the set of first $k$ neighbours of $x$ in $\gamma$. Here, $N(\gamma)$ is the cardinality of $\gamma$. In case there are two or more points within the same distance from a given point, we use the lexicographic ordering as a “tie-breaker” to determine the $k$ nearest neighbor structure. However, such ties have zero probability for the random point sets considered here.

    For a finite configuration $\gamma$, we consider an energy of the form
    $$H(\gamma)= \sum_{x \in \gamma}\sum_{y \in V^k(x,\gamma)}\Phi(x-y),$$
where $\Phi:\R^d\to \R\cup \{\infty\}$ is some measurable function.   
By adding a new point $x$ into the configuration, the nearest neighbour structure will change, but only locally. If $x \in V^k(y,\gamma\cup\{x\})$ for some $y \in \gamma$, then $V^k(y,\gamma\cup\{x\})= V^{k-1}(y,\gamma)\cup\{x\}.$ Otherwise, it is obvious that $V^k(y,\gamma\cup\{x\})=V^k(y,\gamma)$. Due to this fact, the corresponding Papangelou intensity takes the form \eqref{Gibbs} with
\begin{equation}\label{local_energy_NNG}
    h(x,\gamma):=\sum_{y \in V^k(x,\gamma)}\Phi(x-y)+ \sum_{y \in \gamma}\I_{V^k(y,\gamma\cup\{x\})}(x)\left[\Phi(y-x)- \Phi(y-v^{k}(y,\gamma))\right].
\end{equation}

\begin{prop}\label{L8}
There exist a constant $N_d$ depending only on the dimension $d$ such that 
$$z e^{-\beta (1+2 N_d)k\|\Phi\|_{\infty}}\leq \lambda^*(x,\gamma)\leq z e^{\beta (1+2 N_d)k\|\Phi\|_{\infty}}.$$
\end{prop}

\begin{proo}
In \eqref{local_energy_NNG}, the indicator $\I_{V^k(y,\gamma\cup\{x\})}(x)$ takes the value $1$ only for finitely many points $y \in \gamma$. The number of such points is random, yet bounded by $k N_d$, where $N_d$ depends only on the dimension (see Lemma 4.3. in \cite{P03} for the proof). 

\qed
\end{proo}

\begin{coro}\label{C5}
Let $\Phi:\R^d\to \R$ be such that $\|\Phi\|_{\infty}<\infty$. Then the Gibbs point process defined by the local energy \eqref{local_energy_NNG} is non-hyperuniform.
\end{coro}

The Papangelou intensity in Corollary \ref{C5} is double locally stable. Therefore, it is enough to verify the moment condition of Proposition \ref{L1}, (iv). Alternatively, if we assume that $\Phi$ is decreasing and non-negative, we are allowed to have an explosion around the origin. 

\begin{coro}\label{C6}
If $\Phi(x):= \Phi(\|x\|)$ is decreasing and non-negative function on $\R_+$, then the Gibbs point process defined by \eqref{local_energy_NNG} is non-hyperuniform.
\end{coro}

\begin{exam}[Coulomb interaction]
    For $d\geq 3$, let $\Phi(x)= \frac{1}{\|x\|^{d-2}}$. Then by Corollary \ref{C6}, the Gibbs point process defined by \eqref{local_energy_NNG} is non-hyperuniform.
\end{exam}

For the proofs, see Section \ref{S5}.

\section{Proofs of the main results}\label{S5}
This section aims to present the proof of the main result of this paper together with the proofs of our examples and auxiliary lemmas that involve more technical details.

\subsection{Proofs of Theorem \ref{thm1} and Theorem \ref{thm2}}
We start by proving Theorem \ref{thm2}. A slight modification of the proof then yields the statement of Theorem \ref{thm1}.
\begin{proo}[Theorem \ref{thm2}]
 We aim to show that there is a constant $F>-1$ not depending on $\Lambda$ such that
\begin{equation}\label{assert2}
    \frac{\V(N_\Lambda)}{\E N_\Lambda}\geq \frac{1}{1+F}.
\end{equation}
By the stationarity of $\Gamma$, this already implies \eqref{VarBound2}. 

In the spirit of the proof of Proposition $4.1.$ in \cite{R70}, we define $R_\Lambda:= R_\Lambda (\gamma):=  \int_{\Lambda\times\mathbb{M}} \lambda^*(\mathbf{x},\gamma)\nu(\D \mathbf{x})$ for $\gamma \in \mathcal{N}$. One has directly by \eqref{GNZ} that $\E N_\Lambda = \E R_\Lambda$. Moreover, note that
\begin{equation}\label{ERR}
    \E R^2_{\Lambda} =  \int \int_{(\Lambda\times\mathbb{M})^2} \lambda^*(\mathbf{x},\gamma)\lambda^*(\mathbf{y},\gamma)\nu(\D \mathbf{x}) \nu(\D \mathbf{y}) \Gamma(\D \gamma),
\end{equation}
\begin{align}\label{eq5.3}
    \E N_\Lambda (N_\Lambda-1) & = \E \sum_{\mathbf{x} \in \Gamma} \I_{\Lambda\times \mathbb{M}}(\mathbf{x}) \sum_{\mathbf{y} \in \Gamma \setminus \{\mathbf{x}\}} \I_{\Lambda\times\mathbb{M}}(\mathbf{y})\nonumber\\
    & =  \int \int_{(\Lambda\times\mathbb{M})^2} \lambda^*(\mathbf{y},\gamma)\lambda^*(\mathbf{x},\gamma\cup \{\mathbf{y}\})\nu(\D \mathbf{x}) \nu(\D \mathbf{y}) \Gamma(\D \gamma),
\end{align}
where in the latter equality, we used two times \eqref{GNZ}. For the first time, with $f(\mathbf{x},\gamma):= \I_{\Lambda\times\mathbb{M}} (\mathbf{x}) \sum_{\mathbf{y}\in \gamma} \I_{\Lambda\times\mathbb{M}}(\mathbf{y})$ and secondly, with $g(\mathbf{y},\gamma):= \I_{\Lambda\times\mathbb{M}} (\mathbf{y}) \int_{\Lambda\times\mathbb{M}} \lambda^*(\mathbf{x},\gamma\cup \{\mathbf{y}\})\nu(\D \mathbf{x}).$ Similarly, by the definition of $R_\Lambda$ and \eqref{GNZ}, we have
\begin{equation}\label{ENR}
\E N_\Lambda R_\Lambda = \E N_\Lambda (N_\Lambda-1).
\end{equation}
Suppose for a moment that
\begin{equation}\label{assM}
    \mathcal{M}:= \E R^2_\Lambda - \E  N_\Lambda (N_\Lambda-1)\leq F \E N_\Lambda
\end{equation}
for some $F\in (-1,\infty)$, we have that
\begin{align*}
    (F+1)^2 (\E N_\Lambda)^2 & = [\E(F+1)N_\Lambda]^2= [\E(F N_\Lambda + R_\Lambda)]^2 \leq \E(F N_\Lambda + R_\Lambda)^2 \\
    & = F^2 \E N^2_\Lambda + 2 F \E N_\Lambda (N_\Lambda -1) + \E R^2_\Lambda\\
    &= F^2 \E N^2_\Lambda + (2 F+1) \E N_\Lambda (N_\Lambda -1) + \mathcal{M}\\
    &\leq F^2 \E N^2_\Lambda + (2 F+1) \E N_\Lambda (N_\Lambda -1) +F \E N_\Lambda\\
    &= (F+1)^2 \E N_\Lambda^2 - (F+1)\E N_\Lambda.
\end{align*}
Here, we have used consecutively \eqref{GNZ}, Cauchy--Schwartz inequality, \eqref{ENR} and the assumption \eqref{assM}. From this, the assertion \eqref{VarBound2} immediately follows. 

It remains to check the validity of \eqref{assM}.
First, we write $\mathcal{M}=\mathcal{M}_1 + \mathcal{M}_2+ \mathcal{M}_3$, where
\begin{align*}
\mathcal{M}_1 &:=  \int \int_{(\Lambda\times\mathbb{M})^2} \lambda^*(\mathbf{x},\gamma)\lambda^*(\mathbf{y},\gamma)\I\{|x-y|> \delta\}\nu(\D \mathbf{x}) \nu(\D \mathbf{y}) \Gamma(\D \gamma)\\
& - \int \int_{(\Lambda\times\mathbb{M})^2} \lambda^*(\mathbf{y},\gamma)\lambda^*(\mathbf{x},\gamma\cup \{\mathbf{y}\})\I\{|x-y|> \delta\}\nu(\D \mathbf{x}) \nu(\D \mathbf{y}) \Gamma(\D \gamma),\\
\mathcal{M}_2 &:= \int \int_{(\Lambda\times\mathbb{M})^2} \lambda^*(\mathbf{x},\gamma)\lambda^*(\mathbf{y},\gamma)\I\{|x-y|\leq \delta\}\nu(\D \mathbf{x}) \nu(\D \mathbf{y}) \Gamma(\D \gamma),\\
\mathcal{M}_3&:= - \int \int_{(\Lambda\times\mathbb{M})^2} \lambda^*(\mathbf{y},\gamma)\lambda^*(\mathbf{x},\gamma\cup \{\mathbf{y}\})\I\{|x-y|\leq \delta\}\nu(\D \mathbf{x}) \nu(\D \mathbf{y}) \Gamma(\D \gamma).
\end{align*}
Note that
\begin{equation}\label{defM}
\mathcal{M}_1=\E\int_{\Lambda} \int_{\Lambda\setminus B(x,\delta)} \int_{\mathbb{M}^2}\lambda^*(\mathbf{y},\Gamma)\lambda^*(\mathbf{x},\Gamma) \left(1- \frac{\lambda^*(\mathbf{x},\Gamma\cup\{\mathbf{y}\})}{\lambda^*(\mathbf{x},\Gamma)}  \right) \nu(\D \mathbf{x}) \nu(\D \mathbf{y}).
\end{equation}
To proceed further, we denote the quantities 
\begin{align*}
    D_1&:= \int_{\mathbb{M}} \E |\lambda^*((0,m),\Gamma)|^{2\alpha_1}\lambda_{\mathbb{M}}(\D m)\\
    D_2&:=\int_{\R^d\setminus B(0,\delta)}\left(\int_{\mathbb{M}^2}\E \left|1-\frac{\lambda^*((0,m_1),\Gamma \cup \{(y,m_2)\})}{\lambda^*((0,m_1),\Gamma)}\right|^{\alpha_2}\D\lambda_{\mathbb{M}}(m_1,m_2)\right)^{1/\alpha_2} \D y
\end{align*}
that are both assumed to be finite. Using Fubini's theorem, the H{\"o}lder inequality w.r.t. the product measure $P_{\Gamma}\otimes \lambda_{\mathbb{M}}^2$ ($P_{\Gamma}$ is the distribution of $\Gamma)$ and finally, the stationarity of $\Gamma$, we arrive at
\begin{align*}
  &\mathcal{M}_1   =  \int_{\Lambda} \int_{\Lambda\setminus B(x,\delta)} \int_{\mathbb{M}^2} \E \lambda^*(\mathbf{y},\Gamma)\lambda^*(\mathbf{x},\Gamma) \left(1- \frac{\lambda^*(\mathbf{x},\Gamma\cup\{\mathbf{y}\})}{\lambda^*(\mathbf{x},\Gamma)}  \right) \nu(\D \mathbf{x}) \nu(\D \mathbf{y})\\
  & \leq \int_{\Lambda} \int_{\Lambda\setminus B(x,\delta)}\left(\int_{\mathbb{M}^2} \E |\lambda^*((y,m_y),\Gamma)|^{2\alpha_1}\D\lambda_{\mathbb{M}}(m_y,m_x)\int_{\mathbb{M}^2} \E |\lambda^*((x,m_x),\Gamma)|^{2\alpha_1}\D\lambda_{\mathbb{M}}(m_y,m_x)\right)^{\frac{1}{2\alpha_1}}\\
  &\hspace{3cm} \cdot\left(\int_{\mathbb{M}^2}\E\left|1- \frac{\lambda^*((x,m_x),\Gamma\cup\{(y,m_y)\})}{\lambda^*((x,m_x),\Gamma)}  \D\lambda_{\mathbb{M}}(m_x,m_y)\right|^{\alpha_2}\right)^{1/\alpha_2} \D x \D y\\
   & =  \left(\lambda_{\mathbb{M}}(\mathbb{M})D_1\right)^{1/\alpha_1} \int_{\Lambda} \int_{\Lambda\setminus B(x,\delta)}\left(\int_{\mathbb{M}^2}\E\left|1- \frac{\lambda^*((0,m_x),\Gamma\cup\{(y-x,m_y)\})}{\lambda^*((0,m_x),\Gamma)}  \right|^{\alpha_2} \D \lambda^2_{\mathbb{M}}(m_x, m_y)\right)^{1/\alpha_2}\\
   &\hspace{13cm}\times\D x \D y\\
    & \leq \left(\lambda_{\mathbb{M}}(\mathbb{M})D_1\right)^{1/\alpha_1} \int_{\Lambda} \int_{B(0,\delta)^C}\left(\int_{\mathbb{M}^2}\E\left|1- \frac{\lambda^*((0,m_1),\Gamma\cup\{(z,m_2)\})}{\lambda^*((0,m_1),\Gamma)}  \right|^{\alpha_2} \D \lambda^2_{\mathbb{M}}(m_1, m_2)\right)^{1/\alpha_2}\\
   &\hspace{13cm}\times\D x \D z\\  
   & = |\Lambda|\left(\lambda_{\mathbb{M}}(\mathbb{M})D_1\right)^{1/\alpha_1}  D_2 =: F_1' |\Lambda| =: F_1 \E N_\Lambda.
\end{align*}
By Assumptions (A1) and (A2), $F_1 <\infty$.
By the same arguments,
\begin{align*}
    \mathcal{M}_2 & = \int_{\Lambda} \int_{\Lambda\cap B(x,\delta)} \int_{\mathbb{M}^2} \E \lambda^*((x,m_x),\Gamma)\lambda^*((y,m_y),\Gamma)\D \lambda_{\mathbb{M}}^2(m_x,m_y)\D y \D x\\
    & \leq  \int_{\Lambda} \int_{B(x,\delta)} \int_{\mathbb{M}^2} \E |\lambda^*((0,m_x),\Gamma)|^{2}\D \lambda_{\mathbb{M}}^2(m_x,m_y)\D y \D x\\
    & = |\Lambda| |B(0,\delta)| \lambda_{\mathbb{M}}(\mathbb{M}) \int_{\mathbb{M}} \E |\lambda^*((0,m),\Gamma)|^{2} \lambda_{\mathbb{M}}(\D m)=:F_2 \E N_{\Lambda}
\end{align*}
is finite by (A1). At last, $\mathcal{M}_3\leq 0$ and hence \eqref{assM} holds true with $F=F_1 + F_2$.

\qed
\end{proo}

\begin{proo}[Theorem \ref{thm1}]
Clearly, for the constants $\mathcal{M}_2, \mathcal{M}_3$ from the proof of Theorem \ref{thm2}, one can use the same 
 arguments in order to find finite upper bounds. For $\mathcal{M}_1$, we use Cauchy--Schwartz inequality to see that
\begin{align*}
  \mathcal{M}_1  & \leq  \int_{\Lambda} \int_{\Lambda\setminus B(x,\delta)}\phi(x-y) \int_{\mathbb{M}^2} \E \lambda^*(\mathbf{y},\Gamma)\lambda^*(\mathbf{x},\Gamma) \nu(\D \mathbf{x}) \nu(\D \mathbf{y})\\
& \leq \left(\int_{\mathbb{M}^2}\E \lambda^*((0,m_x),\Gamma)^2\D\lambda_{\mathbb{M}}(m_x,m_y)\right)\int_{\Lambda} \int_{\Lambda\setminus B(x,\delta)}\phi(x-y) \D x \D y\\
& \leq \left(\int_{\mathbb{M}^2}\E \lambda^*((0,m_x),\Gamma)^2\D\lambda_{\mathbb{M}}(m_x,m_y)\right)|\Lambda| \int_{B(0,\delta)^C}\phi(z) \D z := F_1 \E N_{\Lambda}.
\end{align*}
The rest of the proof follows the same path as the proof of Theorem \ref{thm2}.

  \qed 
\end{proo}

\subsection{Proof of Proposition \ref{L4}}\label{Sec_P5}
\begin{proo}[]

For all $k\in \Z^d$ we denote by $\mathcal{C}_k$ a cube of a side length $\varepsilon$ lying in the centre of $\mathcal{D}_k$, while both cubes have parallel edges. Now, let $k \in \Z^d$ and $\gamma \in \cN$ be chosen arbitrary, but fixed. Let $\Pi^z_{\mathcal{D}_k}$ be the marked Poisson point process in $\mathcal{D}_k\times \mathbb{M}$ with intensity $z$ and $\cN_{\mathcal{D}_k}\subset \cN$ the set of configurations restricted to $\mathcal{D}_k\times\mathbb{M}$. It can be shown that $\Gamma$ satisfies also DLR equations (see Theorem 1 in \cite{D19}). Then, the distribution of points in $\mathcal{D}_k$ given the configuration in $\mathcal{D}_k^C$ is precisely given by
$$\P(\D \gamma_{\mathcal{D}_k}|\gamma_{\mathcal{D}_k^C})= \frac{1}{Z^{\beta,z}(\gamma_{\mathcal{D}_k^C})}\prod_{n=1}^{N_{\mathcal{D}_k}(\gamma_{\mathcal{D}_k})} \lambda^*(\mathbf{x}_n,\gamma_{\mathcal{D}_k^C}\cup \{\mathbf{x}_1,\ldots, \mathbf{x}_{n-1}\} ) \Pi^z_{\mathcal{D}_k}(\D \gamma_{\mathcal{D}_k}),$$
where
$$Z^{\beta,z}(\gamma_{\mathcal{D}_k^C})= \int_{\cN_{\mathcal{D}_k}} \prod_{n=1}^{N_{\mathcal{D}_k}(\gamma_{\mathcal{D}_k})} \lambda^*(\mathbf{x}_n,\gamma_{\mathcal{D}_k^C}\cup \{\mathbf{x}_1,\ldots, \mathbf{x}_{n-1}\} ) \Pi^z_{\mathcal{D}_k}(\D \gamma_{\mathcal{D}_k}).$$

Since $\lambda^*$ is locally stable from above, 
\begin{align*}
    Z^{\beta,z}(\gamma_{\mathcal{D}_k^C})& \leq \int_{\cN_{\mathcal{D}_k}}{C_1}^{N_{\mathcal{D}_k}(\gamma_{\mathcal{D}_k})} \Pi^z_{\mathcal{D}_k}(\D \gamma_{\mathcal{D}_k})\\
    & = \sum_{n=0}^{\infty}C_1^n\frac{z^n(2\delta+\varepsilon)^{d n}}{n!}e^{-z(2\delta+\varepsilon)^d} = : Z^{\beta,z}<\infty.
\end{align*}
The constant $Z^{\beta,z}$ depends on $\delta, \varepsilon, C_1$ and $d$, but not on $k$ nor the configuration $\gamma_{\mathcal{D}_k^C}$. Then, using that $d(\mathbf{x}, \gamma_{\mathcal{D}_k^C})\geq \delta$ for any $\mathbf{x}\in \mathcal{C}_k\times\mathbb{M}$, we have that

\begin{align*}
\P(&(I^{2\delta+\varepsilon}(\Gamma))_k= 1|\gamma_{\mathcal{D}_k^C}) = \P(N_{\mathcal{D}_k}(\Gamma)\geq 1|\gamma_{\mathcal{D}_k^C}) \\
& \geq \P(N_{\mathcal{D}_k}(\Gamma)=N_{\mathcal{C}_k}(\Gamma)= 1|\gamma_{\mathcal{D}_k^C}) \\
& = \frac{1}{Z^{\beta,z}(\gamma_{\mathcal{D}_k^C})}\int_{\cN_{\mathcal{D}_k}} \I\{N_{\mathcal{D}_k}(\gamma)=N_{\mathcal{C}_k}(\gamma)=1\} \lambda^*(\gamma_{\mathcal{C}_k},\gamma_{\mathcal{D}_k^C})
\Pi^z_{\mathcal{D}_k}(\gamma_{\mathcal{D}_k})\\
& \geq \frac{C_2}{Z^{z,\beta}}  \int_{\cN_{\mathcal{D}_k}}\I\{N_{\mathcal{D}_k}(\gamma)=N_{\mathcal{C}_k}(\gamma)=1\}\Pi^z_{\mathcal{D}_k}(\gamma_{\mathcal{D}_k})\\
&= \frac{C_2}{Z^{z,\beta}} \P(\Pi^z_{\mathcal{D}_k}\cap (\mathcal{C}_k\times\mathbb{M})=1) \P(\Pi^z_{\mathcal{D}_k}\cap(\mathcal{D}_k \setminus \mathcal{C}_k\times\mathbb{M})=0) \\
& = \frac{C_2}{Z^{z,\beta}}z \varepsilon^d e^{-z (2\delta+\varepsilon)^d}=:p >0,
\end{align*}
where $p$ does not depend neither on $k$ nor the boundary condition $\gamma_{\mathcal{D}_k^C}$. Note that the quantity $\lambda^*(\gamma_{\mathcal{C}_k},\gamma_{\mathcal{D}_k^C})$ is well defined, since $\gamma_{\mathcal{C}_k}$ is assumed to be almost surely a one-point set.

To prove the statement, we construct a disagreement coupling of $I^{2\delta+\varepsilon}(\Gamma)$ and $B \sim B(p)^{\otimes \Z^d}$. For this purpose, let $I$ be a finite index set and $X:=(X_i)_{i \in I}$ be random variables, not necessarily independent or identically distributed, with values in $\{0,1\}$. Assume that
\begin{equation}\label{coupling}
\P(X_i=1|X_j; j \in I\setminus \{i\})>p, \quad \forall i \in I.
\end{equation}
Then also
$$\P(X_i=1)= \sum_{x_j \in \{0,1\}; j\in I\setminus \{i\}} \P(X_i=1|X_j=x_j; j \in I\setminus \{i\}) \P(X_j = x_j; j \in I\setminus\{i\})>p$$
and, similarly,
$$\P(X_i=1|X_j; j\in J)>p$$
for any $J\subseteq I\setminus\{i\}$.

Let $U:=(U_i)_{i \in I}$ be a vector of i.i.d  uniform random variables on $[0,1]$ and define $B_{i}:= \I[U_{i}<p], i \in I$. Clearly, $B_I:=(B_i)_{i \in I}$ is a vector of i.i.d Bernoulli variables with parameter $p$. 

For the coupling, we define $Z_{i_1}= \I[U_{i_1} < \P(X_{i_1}=1)], i_1 \in I$. It is easy to check that $\P(Z_{i_1}= z)=\P(X_{i_1}= z)$ for $z \in \{0,1\}$ and $\P(B_{i_1}\leq Z_{i_1})=1$. 
Inductively for $k\geq 1$, let $((Z_{i_1},\ldots Z_{i_k}), (B_{i_1},\ldots B_{i_k}))$ be a coupling of $(X_{i_1},\ldots X_{i_k})$ and $(B_{i_1},\ldots B_{i_k})$, such that $\P(Z_{i_j}\geq B_{i_j})=1$ for any $j=1,\ldots,k$. Then we define
$$Z_{i_{k+1}}:= \I[U_{i_{k+1}}<\P(X_{i_{k+1}}|X_{i_j}=Z_{i_j}; j=1,\ldots,k)].$$
Again, $\P(Z_{i_{k+1}}\geq B_{i_{k+1}})=1$ and for $z_1,\ldots z_k \in \{0,1\}$, we compute
\begin{align*}
        \P(Z_{i_{1}}=z_1,&\ldots, Z_{i_k}=z_k, Z_{i_{k+1}}=1)\\
    &= \P(Z_{i_{1}}=z_1,\ldots, Z_{i_k}=z_k)\P(Z_{i_{k+1}}=1|Z_{i_{1}}=z_1,\ldots, Z_{i_k}=z_k)\\
    & = \P(X_{i_{1}}=z_1,\ldots, X_{i_k}=z_k)\P(U_{i_{k+1}}<\P(X_{i_{k+1}}=1|X_{i_{1}}=z_1,\ldots, X_{i_k}=z_k))\\
    &= \P(X_{i_{1}}=z_1,\ldots, X_{i_k}=z_k)\P(X_{i_{k+1}}=1|X_{i_{1}}=z_1,\ldots, X_{i_k}=z_k)\\
    &= \P(X_{i_{1}}=z_1,\ldots, X_{i_k}=z_k, X_{i_{k+1}}=1).
\end{align*}
Let $Z:= (Z_i)_{i \in I}$. We showed that $(Z,B_I)$ is a coupling of $X$ and $B_I$ yielding
$$B_I< < X.$$
It remains true if we turn to the limit and take $I_n \nearrow \Z^d$. The choice $X = (I^{2\delta+\varepsilon}(\Gamma))_{i \in \Z^d}$ satisfies the assumption \eqref{coupling} and hence, by the arguments above
$$B< < I^{2\delta+\varepsilon}(\Gamma).$$

\qed
\end{proo}

 \subsection{Proof of Corollary \ref{C4}}
\begin{proo}[]
As before, one needs to verify Assumptions (A1), (A2) of Theorem \ref{thm2}. Here, it is enough to do so in the form without marks.

By Proposition \ref{L7}, the Papangelou intensity of $\Gamma$ is locally stable from below and hence, by Proposition \ref{L_PoissonSandwich}, there exists a Poisson point process $\Pi$ such that $\Pi \ll \Gamma$ and $\Pi$ has the intensity $z e^{-K \beta}$, where $K$ is 
 as in \eqref{localEnergyBound2} and $z, \beta$ are the parameters defining the Papangelou intensity \eqref{Gibbs} of the Gibbs point process $\Gamma$. This stochastic minoration can be interpreted via coupling of $\Pi$ and $\Gamma$ such that $\gamma'\subseteq \gamma$ whenever $\gamma'\sim\Pi$ and $\gamma \sim \Gamma$. For Voronoi tessellation, it translates to $C(x,\gamma)\subseteq C(x,\gamma')$ for any $x \in \gamma\cap\gamma'$. Together with the upper bound from Proposition \ref{L7}, we may write
 \begin{align*}
     \E|\lambda^*(0,\Gamma)|^{2\alpha_1} & \leq z^{2\alpha_1} e^{2\alpha_1 \beta K} \E e^{2\alpha_1 \beta |C(0,\Pi)|}.
 \end{align*}

 Now we aim to construct a ball $B(0,R)$ such that $R=R(\Pi)$ is a random variable depending on $\Pi$ and $C(0,\Pi)\subseteq B(0,R)$ a.s. We use the exact approach as in the proof of Lemma 5.1 in \cite{P07} with the only difference that we are interested in more precise estimates of the tail probabilities $\P(R>t), t>0$.

Let $K_1,\ldots, K_J$ be set of circular cones with apices in the origin with angular radii $\pi/6$. We do not expect the cones to have zero volume intersections, yet we chose $J$ to be minimum such that $\cup_{j=1}^J K_j =\R^d$. Note that $J$ depends on the dimension and is finite. For each $j=1,\ldots,J$ chose $x_j \in \gamma' \cup K_j$ to be the closest point to the origin. 

Denote by $H_x(y):= \{z \in \R^d; \|z-y\|\leq \|z-x\|\}$ the closed half-space induced by points that are closer to $y$ than to $x$. By definition of Voronoi cell,
$$C(0,\Pi)=\bigcap_{x \in \gamma'}H_x(0)\subseteq \bigcap_{j=1}^J H_{x_j}(0)\subseteq \bigcup_{j=1}^J H_{x_j}(0) \cap K_j \quad\text{a.s}.$$
Set $R:=\max\{\|x_j\|; j=1,\ldots,J\}$. We need to verify that $C(0,\Pi) \subseteq B(0,R)$ a.s. To do so, we chose $y \in H_{x_j}(0) \cap K_j$ and show that $\|y\|\leq \|x_j\|$ for any $j\in \{1,\ldots,J\}$. It follows from a simple computation
$$\|y\|^2\leq \|y-x_j\|^2=\|x_j\|^2+\|y\|^2 -2\langle x_j,y \rangle \leq \|x_j\|^2+\|y\|^2 - \|x_j\| \|y\|,$$
where the first inequality holds since $y \in H_{x_j}(0)$ and the last one from the fact that $y \in K_j$ (note that $z_1,z_2 \in K_j$ implies $\langle z_1, z_2\rangle \geq 1/2\|z_1\|\|z_2\|$). Finally, we have that $C(0,\Gamma)\subseteq B(0,R(\Pi))$ a.s.

Using void probabilities of the Poisson point process $\Pi$, we arrive at the estimate 
\begin{align*}
    \P(R>t)& \leq \sum_{j=1}^J \P(\|x_j\|>t)=\sum_{j=1}^J\P(\Pi \cap B(0,t)\cap K_j = \emptyset)\\
    & = \sum_{j=1}^J \exp\{-z e^{-\beta K}|K_j\cap B(0,t)|\}\\
    & = J \exp\{-z e^{-\beta K}t^dc_d\},
\end{align*} 
where $c_d := |K_1\cap B(0,1)|_d$. We let the reader check that 
$$c_d = |B(0,1)|_{d-1}\left(\sin^{d-1}(\pi/12)\cos(\pi/12)\frac{1}{d}+\int_{0}^{\pi/12}\sin^d(\theta)\D \theta\right).$$

Eventually, for (A1), we have that 
\begin{align}\label{A1vor}
    \E\lambda^*(0,\Gamma)^{2\alpha_1}&\leq z^{2\alpha_1} e^{2\alpha_1 \beta K}\E_{\Pi}e^{2\alpha_1\beta |B(0,1)|_d R^d}\\
    & = z^{2\alpha_1} e^{2\alpha_1 \beta K} \int_{0}^{\infty} \P(\exp\{2\alpha_1\beta |B(0,1)|_d R^d\}>t) \D t\nonumber\\
   & = z^{2\alpha_1} e^{2\alpha_1 \beta K} d\int_{0}^{\infty} u^{d-1} \exp\{2\alpha_1\beta |B(0,1)|_d u^d\} \P(R>u) \D u\nonumber\\
   &\leq  z^{2\alpha_1} e^{2\alpha_1 \beta K} d\int_{0}^{\infty} u^{d-1} \exp\{2\alpha_1\beta |B(0,1)|_d u^d\}  J \exp\{-z e^{-\beta K}u^dc_d\} \D u\nonumber.
\end{align}
The latter integral converges as long as
\begin{equation}\label{alpha1}
\beta < \frac{1}{2\alpha_1 |B(0,1)|_d} z e^{-\beta K} c_d.
\end{equation}

Similarly for (A2), we have from the assumptions on $\Phi$ that
$$0\leq \frac{\lambda^*(0,\gamma\cup\{y\})}{\lambda^*(0,\gamma)}\leq e^{2\beta K} e^{\beta|C(0,\gamma)|}, \quad \Gamma-\text{a.s.}$$
Therefore, by an additional use of Fubini's theorem and $1/\alpha_2\leq 1$,

\begin{align*}
     \int_{\R^d}&\left(\E \left|1-\frac{\lambda^*(0,\Gamma \cup \{y\})}{\lambda^*(0,\Gamma)}\right|^{\alpha_2}\right)^{1/\alpha_2} \D {y}\\
    & \leq \int_{\R^d}\E \left(\frac{\lambda^*(0,\Gamma \cup \{y\})}{\lambda^*(0,\Gamma)}\right)^{\alpha_2} \I\{|y| \leq R \}\D {y}\\
    &\leq  \E \int_{B(0,R)}e^{2\alpha_2\beta K} e^{\alpha_2 \beta|B(0,R)|_d} \D y\\   
    &= e^{2\alpha_2\beta K}|B(0,1)|_d\E  R^d e^{\alpha_2 \beta|B(0,1)|_d R^d}.
\end{align*}
A computation in the same spirit as in \eqref{A1vor} shows that the expectation above is finite as long as
\begin{equation}\label{alpha2}
 \beta<\frac{1}{\alpha_2|B(0,1)|_d}z e^{-\beta K}c_d.   
\end{equation}
The two inequalities \eqref{alpha1} \eqref{alpha2} are optimal for $\alpha_1=\frac{3}{2}$ and $\alpha_2=3$.

\qed
\end{proo}

    \begin{rema}
     We shall highlight two facts for the construction of the tail probabilities estimates of $R$.
       First, the construction of the estimates is not damaged by the fact that the cones $K_1,\ldots, K_J$ may intersect and we can possibly choose some point $x \in \gamma'$ multiple times. Second, it is seen from the exponential form of the tail probabilities that the circular cones are the optimal choice for this construction. They give us more precise estimates than any other solids would do (e.g. non circular and non intersecting). 
    \end{rema}

\subsection{Proofs of Corollary \ref{C5} and Corollary \ref{C6}}
\begin{proo}[Corollary \ref{C5}]
The Papangelou intensity is double locally stable by Proposition \ref{L8}, ergo Assumption (A1) of Theorem \ref{thm2} is justified by Proposition \ref{L1}, (i).

It remains to validate (A2). For that, we construct $R:= R(\gamma)$ such that $R$ is the range of interaction, i.e.
$$\lambda^*(0,\gamma)=\lambda^*(0,\gamma \cap B(0,R)) \quad \Gamma-\text{a.s.}$$
The construction is almost the same as in the proof of Lemma 6.1. in \cite{PY01}, where the authors proved stabilization for a planar undirected $k$-nearest neighbour graph (there is an edge between $x,y$ whenever $x \in V^k(y,\gamma)$ or $y \in V^k(x,\gamma)$). Here, we consider a directed graph, i.e. such that there is an edge pointing from $x$ to $y$ whenever $y\in V^k(x,\gamma)$. The idea, however, remains the same. 

Let $J$ be the smallest integer such that $K_1,\ldots,K_J$ are cones with the apex at the origin and the angular radius at most $\pi/6$ such that $\cup_{j=1}^J K_j= \R^d$. Note that, compared to the proof of Lemma \ref{C4}, we do not need to optimize the shape of the cones. In fact, by Proposition \ref{L1}, (iv), all we need to show is that the radius of has finite $\alpha$ moment for some $\alpha>d$ and set $\alpha_2=\alpha/d$. Let $D:=D(\gamma)$ be the smallest $t>0$ such that there are at least $k+1$ points in each $K_j\cap B(0,t), j=1,\ldots,J$ and set $R=2D$.

First, for any $y\in \gamma$ such that $0\in V^k(y,\gamma)$ there exists $j\in\{1,\ldots,J\}$ such that $y \in K_j\cap B(0,D)$. Otherwise, $y$ would have at least $k+1$ points that are closer to $y$ than the origin by the construction of $D$ and that contradicts $0\in V^k(y,\gamma)$. In addition, $V^k(y,\gamma)\in B(0,R)$. No point from $B(0,R)^C$ can be among the $k$ nearest neighbours of $y$, because $y\in K_j\cap B(0,D)$ implies that there are at least $k$ points closer to $y$ than a potential neighbour outside $B(0,R)$. We conclude that $R$ is a decreasing range of interaction for the Papangelou intensity $\lambda^*$ as in Definition \ref{Def_range}.

By Proposition \ref{L8} and Proposition \ref{L_PoissonSandwich}, there is a Poisson point process $\Pi$ with intensity $\lambda:=z e^{-\beta (1+2 N_d) k N_d\|\Phi\|_{\infty}}$ such that $\Pi \ll \Gamma$. By Corollary \ref{C01}, for any $t>0$,

\begin{align*}
    \P(D(\Gamma)>t)&\leq \P(D(\Pi)>t)\leq \sum_{j=1}^J \P(\#(\Pi\cap T_j \cap B(0,t))\leq k)\\
    &= \sum_{j=1}^J \sum_{i=1}^k \frac{(\lambda|T_j\cap B(0,t)|)^i}{i!} e^{-\lambda |T_j\cap B(0,1)| t^d}
\end{align*}
Let $\alpha>d$. Then,
\begin{align*}
    \E R(\Gamma)^\alpha &\leq  \E R(\Pi)^\alpha = \int_{0}^{\infty}\P(R(\Pi)^\alpha>t)\D t\\
    &= \int_{0}^{\infty}\P((2 D(\Pi))^\alpha>t)\D t \\
    & = \int_{0}^{\infty}2^\alpha u^{\alpha-1}\P(D(\Pi)>u)\D u \\
    &\leq \sum_{j=1}^J \sum_{i=1}^k \frac{2^\alpha}{i!}\int_0^{\infty} u^{\alpha-1} (\lambda|T_j\cap B(0,u)|)^i e^{-\lambda |T_j\cap B(0,1)| u^d}\D u.
\end{align*}
The latter term is finite as it is a finite sum of converging integrals. The Assumption (A2) of Theorem \ref{thm2} is satisfied for $\alpha_2=\alpha/d$ by Proposition \ref{L1}, (iv).  

\qed
\end{proo}

\begin{proo}[Corollary \ref{C6}]
Since $\Phi\geq 0$ everywhere, the Papangelou intensity is local stable from above, hence Assumption (A1) of Theorem \ref{thm2} is trivially satisfied by Proposition \ref{L1}, (i). For (A2), an easy computation leads to
\begin{align*}
    \frac{\lambda^*(0,\gamma \cup \{y\})}{\lambda^*(0,\gamma)} = & \exp \left\lbrace  -\beta\left[\Phi(y)-\Phi(v^{k+1}(0,\gamma \cup\{0,y\}))\right]\I\{y \in V^k(y,\gamma \cup\{0,y\})\} \right.\\
    & \hspace{1cm} -\beta \left[\Phi(y)-\Phi(y-v^{k+1}(y,\gamma \cup \{0,y\}))\right]\I\{0\in V^k(y,\gamma \cup \{0,y\})\}\\
    &\hspace{1cm} +\beta \sum_{\substack{x \in \gamma \\ y \in V^k(x,\gamma \cup \{0,y\})\\ 0 = v^{k+1}(x,\gamma \cup\{0,y\})}} \left.\left[\Phi(x)-\Phi(x-v^{k+2}(x,\gamma \cup \{0,y\}))\right]\right\rbrace\\
    &\leq \exp\left\lbrace +\beta \sum_{\substack{x \in \gamma \\ y \in V^k(x,\gamma \cup \{0,y\})\\ 0 = v^{k+1}(x,\gamma \cup\{0,y\})}} \Phi(x) \right\rbrace,
\end{align*}
where we used the fact that $\Phi(x)$ is non-negative and decreasing. Let $\delta>0$ and $\zeta <\infty$ be such that $\Phi(x)<\zeta$ whenever $\|x\|\geq \delta/2$. Then, for $\alpha_2>1$,
\begin{align*}
   D &:=\int_{\R^d\setminus B(0,\delta)}\left(\E \left|1-\frac{\lambda^*(0,\Gamma' \cup \{y\})}{\lambda^*(0,\Gamma')}\right|^{\alpha_2}\right)^{1/\alpha_2} \D {y}\\
    & \leq \int_{\R^d\setminus B(0, \delta)}\left[\E\left(\exp\left\lbrace \beta \sum_{\substack{x \in \gamma \\ y \in V^k(x,\Gamma \cup \{0,y\})\\ 0 = v^{k+1}(x,\Gamma \cup\{0,y\})}} \Phi(x) \right\rbrace - 1\right)^{\alpha_2}\right]^{1/\alpha_2}\D y \\
    & \leq \int_{\R^d\setminus B(0, \delta)}\left[\E\exp\left\lbrace \beta \alpha_2 N_d \zeta \right\rbrace \right.\\
    &\hspace{3cm} \cdot\left.\I\{\exists x \in \Gamma; y \in V^k(x,\Gamma \cup \{0,y\}), 0 = v^{k+1}(x,\Gamma \cup\{0,y\})\}\right]^{1/\alpha_2}\D y\\
    &= A  \int_{\R^d\setminus B(0, \delta)} \left(\P(\exists x \in \Gamma; y \in V^k(x,\Gamma \cup \{0,y\}), 0 = v^{k+1}(x,\Gamma \cup\{0,y\}))\right)^{1/\alpha_2} \D y,
\end{align*}
 where $A:= \exp\{\beta N_d \zeta \}<\infty$. In the third line of the latter expression, we used that $\|x\|\geq \|y\|/2 \geq \delta/2$. Otherwise, $\|x-0\|\leq \|x-y\|$ implies that if $0 = v^{k+1}(x,\Gamma \cup\{0,y\})$, then $y$ cannot be among the $k$ nearest neighbours of $x$.

Again, we let $J$ be the smallest integer such that $K_1,\ldots,K_J$ are cones with the apex at the origin and angular radius at most $\pi/6$ such that $\cup_{j=1}^J K_j= \R^d$. For $\gamma \in \cN$, we define $D:=D(\gamma)$ as the smallest $t>0$ such that there is at least $k+1$ points of $\gamma$ in each $K_j\cap B(0,t), j=1,\ldots,J$ and set $R=2D$. Then 
$$\P(\exists x \in \Gamma; y \in V^k(x,\Gamma \cup \{0,y\}), 0 = v^{k+1}(x,\Gamma \cup\{0,y\})) \leq \P(\|y\|\leq R).$$
If $\|y\|> R$, then any point $x \in \gamma$ with $y \in V^k(x,\gamma\cup \{0,y\})$ and $0=v^{k+1}(x, \gamma\cup \{0,y\})$ satisfies $x \in B(0, D)^C$. But then there is $j \in \{1,\ldots, J\}$ with $x \in K_j$ and at least $k+1$ points in the cone $K_j$ that are closer to $x$ then the origin, contradicting $0 = v^{k+1}(x,\gamma\cup\{0,y\})$.

We are in a position to apply Corollary \ref{C02}. It can be seen from \eqref{local_energy_NNG}, that $\lambda^*(x,\gamma)\leq z$. This is due to the fact that $\Phi$ is non-negative and decreasing. Furthermore, if $d(x,\gamma)\geq \delta/2$, then $\lambda^*(x,\gamma)\geq z \exp\{- (1+N_d)k \zeta\}$. Here, we also used the fact that the number of summands is bounded by a deterministic constant (see Proposition \ref{L8}). By Proposition \ref{L4}, for every $\varepsilon>0$, there is $p>0$ such that $I^{2\delta+\varepsilon}(\Gamma)$ is minorated by a Bernoulli field $B$ with parameter $p$.

For $l \in \mathbb{L}^d:= \{0,1\}^{\Z^d}$, we define $R'(l)$ by taking the smallest $t>0$ such that for all $j=1,\ldots,J$ the cone $K_j$ contains fully at least $k+1$ cubes $\mathcal{D}_{i_1},\ldots, \mathcal{D}_{i_{k+1}}\in \mathcal{D}$ such that $l_{i_1},\ldots, l_{i_{k+1}}=1$. Then $R'$ is decreasing.
By the construction of $I^{2\delta+\varepsilon}(\Gamma)$ and Proposition \ref{L4} we have that almost surely
$$R(\Gamma)\leq R'(I^{2\delta+\varepsilon}(\Gamma)) \leq R'(B)$$
and hence, by Corollary \ref{C02}, 
$$\P(R(\Gamma) > r)\leq \P(R'(I^{2\delta+\varepsilon}(\Gamma)) > r)\leq \P(R'(B)> r).$$

 For $j \in \{1,\ldots,J\}$ take the closest cube $\mathcal{D}_{i_1}\subset K_j$ to the origin and denote by $c_j:= \inf_{x \in \mathcal{D}_{i_1}} d(x,0)$ its distance to the origin. Then there exists another cube $\mathcal{D}_{i_2}\subset K_j$ sharing exactly one vertex with $\mathcal{D}_{i_1}$ in a distance $c_1+c_2$ from the origin, where $c_2$ is the body diagonal length of $\mathcal{D}_{i_1}$. Note that $c^j_1$ and $c_2$ depend on $d, \delta$ and $\varepsilon$ and $c^j_1$ moreover on $j$. Inductively, we construct a chain of cubes $\mathcal{D}_{i_j}; j \in \N$ all fully included in the cone $K_j$. Let $B=(B_i)_{i \in \Z^d}$ be distributed according to $B(p)^{\otimes \Z^d}$. Define $R^j_{diag}(B)$ as the smallest $t= c^j_1+c_2 q$ such that $\sum_{j=1}^q B_{i_j}=k+1$. To simplify the notation, for the rest of the proof let $C$ be a universal finite constant depending on $d, k, \delta, \varepsilon, \alpha_2, c_2$ and $c^j_1,j=1,\ldots,J$.  It can be seen that,

\begin{align*}
    \P( R'(B)>r)&\leq \sum_{j=1}^J \P(R^j_{diag}(B)>r)\\
 & = \sum_{j=1}^J \P(\sum_{j=1}^{\lceil\frac{r-c^j_1}{c_2}\rceil} B_{i_j}\leq k)\\
 & = \sum_{j=1}^J \sum_{n=1} ^k {\binom{\lceil\frac{r-c^j_1}{c_2}\rceil}{n}}p^n (1-p)^{\lceil\frac{r-c^j_1}{c_2}\rceil-n}\\
 &\leq C r^k (1-p)^{r}.
\end{align*}

Eventually, since $p \in (0,1)$,
\begin{align*}
    D &\leq C \int_{\R^d\setminus B(0,\delta)} (\P(R'(B)\geq\|y\|))^{1/\alpha_2}\D y\\
    &\leq C \int_{\R^d} \|y\|^{(k/\alpha_2)}((1-p)^{1/\alpha_2})^{\|y\|}\D y <\infty.
\end{align*}

\qed
\end{proo}

\section*{Acknowledgement}
This work was supported in part by the Labex CEMPI (ANR-11-LABX-0007-01), the ANR project RANDOM (ANR-19-CE24-0014), the CNRS GdR 3477 GeoSto and by the Czech Science Foundation, project no. 22-15763S.

The authors would like to thank the anonymous Referees and the Editor handling our paper for their careful reading and comments that improved the paper.

\end{document}